\documentclass[12pt,a4paper
]{article}

\usepackage[mathscr]{eucal}
\usepackage{enumerate}
\usepackage{epsfig,latexsym,amsfonts,amssymb,amsmath,amscd,
graphics,
epic}
\newcommand{\bp}{\boldsymbol{p}}\newcommand{\bd}{\boldsymbol{d}}\newcommand{\bz}{\boldsymbol{z}}\newcommand{\bn}{\boldsymbol{n}}
\newcommand{\N}{\ensuremath{\mathbb{N}}}
\newcommand{\Z}{\ensuremath{\mathbb{Z}}}

\newcommand{\R}{\ensuremath{\mathbb{R}}}
\newcommand{\C}{\ensuremath{\mathbb{C}}}

\newcommand{\PP}{\ensuremath{\mathbb{P}}}
\newcommand{\redu}
{\ensuremath{\mathrm{red}\:}}\newcommand{\conv}
{\ensuremath{\mathrm{conv}}}\newcommand{\cheb}
{\ensuremath{\mathrm{cheb}}}\newcommand{\jet}
{\ensuremath{\mathrm{jet}}}\newcommand{\Sig}{{\Sigma}}\newcommand{\alp}{{\alpha}}\newcommand{\bet}{{\beta}}
\newcommand{\Int}{\ensuremath{\mathrm{Int}}}

\newcommand{\Idim}{\ensuremath{\mathrm{idim}}}
\newcommand{\Sing}{\ensuremath{\mathrm{Sing}}}\newcommand{\eps}{\ensuremath{\varepsilon}}
\newcommand{\Pic}{\ensuremath{\mathrm{Pic}}}

\newcommand{\lcm}{\ensuremath{\mathrm{lcm}}}

\newcommand{\pr}{\ensuremath{\mathrm{pr}}}\newcommand{\del}{{\delta}}
\newcommand{\gam}{{\gamma}}

\newcommand{\proofend}{\hfill$\Box$\bigskip}

\newtheorem{theorem}{Theorem}[section]
\newtheorem{lemma}{Lemma}[section]
\newtheorem{proposition}[lemma]{Proposition}
\newtheorem{corollary}[lemma]{Corollary}
\newtheorem{definition}[lemma]{Definition}
\newtheorem{remark}[lemma]{Remark}
\newtheorem{example}[lemma]{Example}

\font\ncsc=cmcsc10  \font\ntt=cmtt12
\begin{document}
\title{Counting curves of any genus on $\PP^2_6$}

\author{Mendy Shoval\and Eugenii Shustin}

\maketitle

\begin{abstract}
We obtain a formula for the degrees of the varieties parameterizing
complex algebraic curves of any divisor class and genus on
$\PP^2_6$, the projective plane blown-up at $6$ generic points.
Moreover, the formula computes the degrees of the varieties
parameterizing curves on $\PP^2_6$ which additionally satisfy
certain tangency conditions to a fixed exceptional divisor
$E\subset\PP^2_6$. Our formula contains as special cases the degrees
of the analogous varieties parameterizing curves on $\PP^2_q$, for
$0 \leq q \leq 5$, and on $(\PP^1)^2$. It is an extension of the
Caporaso-Harris formula \cite{CH} counting curves of any degree and
genus in the projective plane, and it can be viewed as an
alternative to the Vakil formula \cite{Va}.
\end{abstract}

\tableofcontents

\section{Introduction}

The celebrated Kontsevich formula \cite{KM} computes recursively all
genus zero Gromov-Witten invariants of the plane. In 1998 L.
Caporaso and J. Harris \cite{CH} suggested a formula recursively
computing the Gromov-Witten invariants of the projective plane for
any genus. These invariants can be viewed as the numbers of plane
curves of degree $d$ and genus $g$ passing through $3d+g-1$ generic
points in the plane, or, equivalently, as degrees of the Severi
varieties in $|{\mathcal O}_{\PP^2}(d)|$ parameterizing plane curves
of degree $d$ and genus $g$. Vakil \cite{Va} extended the
Caporaso-Harris approach further, producing a formula which computes
the number of curves of a given degree and genus in the plane or on
a ruled surface and, in addition, has few fixed multiple points. In
particular, Vakil computes the Gromow-Witten invariants of the plane
blown up at $q\le5$ generic points and the numbers of curves in a
given linear system and of a given genus in the plane blown up at
$6$ points on a conic. In the latter case, using results of Graber
in Gromov-Witten theory \cite{G}, Vakil computes the invariants of
$\PP^2_6$ indirectly, by computing invariants of a degeneration.

The goal of the paper is to prove a direct Caporaso-Harris type
recursive formula for the Gromov-Witten invariants of $\PP^2_6$. The
same formula, in fact, computes the Gromov-Witten invariants of all
Del Pezzo surfaces of degree $\ge3$. Besides a direct enumerative
outcome, we had in mind another application: our formula has a real
version, which computes (purely real) Welschinger invariants of all
real Del Pezzo surfaces of degree $\ge3$ \cite{IKS,IKS1},
particularly, of the two-component real cubic surface in $\PP^3$ (so
far we did not succeed to do such computation when using the pair
(plane, conic) as in \cite[Section 9.2]{Va}).

We follow the reasoning of \cite{CH} and \cite{Va} which up to some
extent smoothly applies to the case of $\PP^2_6$, notably, it
allowed us to prove the key bound for the dimension of the moduli
space of curves of a given genus and divisor class, matching given
tangency conditions to a given $(-1)$-curve $E$, and that, for
almost all divisor classes, the general element of such a moduli
space is an immersed (even nodal) curve. The geometric background
for the recursive formula is as follows:
\begin{itemize}\item one computes the degree of a moduli space $V$ in question by
counting curves $C\in V$ passing through $n=\dim V$ generic points
outside $E$, \item specializing one of these points to $E$, one
obtains certain degenerate curves, and, finally, \item computing how
many curves $C\in V$ appear in a deformation of degenerate curves,
one obtains the coefficients in the desired recursive formula.
\end{itemize} The degeneration step again can be done using the
ideas of \cite{CH} (cf. Proposition \ref{pdeg1}), and the result is
that either the degenerate curve remains irreducible while one of
its moving tangency points with $E$ becomes fixed, or the degenerate
curve splits off one copy of $E$ and a number of components
belonging to similar moduli spaces of lower dimension. We should
like to mention that, for $\PP^2_q$, $q\le5$, the reducible
degenerate curves remain \emph{reduced} (cf. \cite[Lemma 5]{IKS}),
and thus, the deformation argument like in \cite{CH} (cf.
\cite[Section 6.2]{Va}) gives the required coefficients
(multiplicities). In turn, for $\PP^2_6$, one encounters a \emph{new
phenomenon}: the reducible degenerate curves can be
\emph{non-reduced}. This requires an additional deformation theory
argument based on the local tropical geometry (section \ref{sec5}),
and this is the main novelty of the present paper.

Finally, we notice that a combination of the approach of the current
paper with the ideas of Vakil to put six blown-up points on a conic
leads to a Capraso-Harris type recursive formula for $\PP^2_7$. It
will be presented in a forthcoming paper.

\textbf{Acknowledgments.} The first author thanks Ya. Varshavsky for
useful conversations and remarks and the Hebrew University for
providing the necessary support and warm working atmosphere. The
second author has enjoyed a support from the Israeli Science
Foundation grant no. 448/09 and from the Hermann-Minkowski-Minerva
Center for Geometry at Tel Aviv University. We also thank I.
Itenberg and V. Kharlamov for very valuable discussions, and R.
Vakil for important remarks.

\section{Families of curves on $\PP^2_6$}

Here we intend to establish several properties of families of curves
on $\PP^2_6$ which will appear in the count, notably, un upper bound
to the dimension, properties of their generic elements, and the
description of zero-dimensional families.

\subsection{General setting}\label{sec-gs}

Let $\Sig=\PP^2_6$ be the complex projective plane blown up at $6$
generic points, $E_i$, $1\le i\le 6$, the exceptional curves of the
blow up, $L\in\Pic(\Sig)$ the pull-back of a generic line in
$\PP^2$, and $E$ a smooth rational curve linearly equivalent to
$2L-E_1-...-E_5$ (the strict transform of a plane conic through the
first five blown up points). Denote by $\Pic_+(\Sig,E)$ and
$\Pic(\Sig,E)$ the semigroups generated by the classes
$D\in\Pic(\Sig)$ of effective irreducible divisors such that $DE> 0$
or $DE\ge 0$, respectively. We have $K_\Sig=-3L+E_1+...+E_6$.

Following \cite[Section 2]{Va}, for a given effective divisor class
$D\in\Pic(\Sig)$ and nonnegative integers $g,n$, denote by
$\overline{\mathcal M}_{g,n}(\Sig,D)$ the moduli space of stable
maps $\bn:\hat C\to\Sig$ of $n$-pointed, connected curves $\hat C$
of genus $g$ such that $\bn_*\hat C\in|D|$ \footnote{In what
follows, by $\bn_*\hat C$ we denote the $\bn$-image of $\hat C$,
whose components are taken with respective multiplicities, by
$\bn(\hat C)$ we denote the reduced image.}. This is an algebraic
stack, whose open dense subset is formed by elements with a smooth
curve $\hat C$ of genus $g$, and the other elements have $\hat C$ of
arithmetic genus $g$ with at most nodes as singularities.

Denote by $\Z^\infty_+$ the direct sum of countably many additive
semigroups \mbox{$\Z_+=\{m\in\Z\ |\ m\ge 0\}$}, labeled by the
naturals, and denote by $e_k\in\Z^\infty_+$ the generator of the
$k$-th summand, $k\in\N$. For
$\alp=(\alp_1,\alp_2,...)\in\Z^\infty_+$, put
$$
\|\alp\|=\sum_{k\ge1}\alp_k,\quad I\alp=\sum_{k\ge1} k\alp_k,\quad
I^\alp=\prod_{k\ge1} k^{\alp_k},\quad\alp!=\prod_{k\ge1}\alp_k!\ ,
$$ and for
$\alp,\alp^{(1)},\ldots,\alp^{(m)} \in \Z_+^\infty$ such that $\alp
\geq \alp^{(1)} + \ldots + \alp^{(m)}$ we put
$$\binom{\alp}{\alp^{(1)},\ldots,\alp^{(s)}} =
\frac{\alp!}{\alp^{(1)}! \ldots \alp^{(m)}!(\alp - \alp^{(1)} - ...
- \alp^{(m)})!}\ .$$ Let a divisor class $D\in \Pic(\Sig,E)$, an
integer $g$, and two elements $\alp,\bet\in\Z^\infty_+$ satisfy
\begin{equation}
0\le g\le g(\Sig,D)=\frac{D^2+DK_\Sig}{2}+1,\quad I\alp+I\bet=DE\ .
\label{en1}
\end{equation}
Given a sequence $\bp=(p_{ij})_{i\ge1,1\le j\le\alp_i}$ of
$\|\alp\|$ distinct points of $E$, define ${\mathcal V}_\Sig(D, g,
\alp, \bet, \bp)\subset\overline{\mathcal M}_{g,\|\alp\|}(\Sig,D)$
as the (stacky) closure of the set of elements $\{\bn:\hat
C\to\Sig,\hat\bp\}$, $\hat\bp=(\hat p_{ij})_{i\ge1,1\le
j\le\alp_j}$, subject to the following restrictions:
\begin{itemize}
\item $\hat C$ is smooth, $\bn(\hat p_{ij})=p_{ij}$ for all $i\ge1$, $1\le j\le\alp_i$,
\item $\bn^*(C\cap E)$ is the following divisor on $\hat C$: \begin{equation}\bn^*(C\cap E)=
\sum_{i\ge1,\ 1\le j\le\alp_i}i\cdot\hat p_{ij}+ \sum_{i\ge1,\ 1\le
j\le\bet_i}i\cdot\hat q_{ij}\ ,\label{ediv}\end{equation} where
$\{\hat q_{ij}\}_{i\ge1,1\le j\le\bet_i}$ is a sequence of
$\|\bet\|$ distinct points of $\hat C\setminus\hat\bp$.
\end{itemize} Put
$$R_\Sig(D,g,\bet)=-D(E+K_\Sig)+\|\bet\|+g-1\ .$$ At last, for a
component $V$ of ${\mathcal V}_\Sig(D,g,\alp,\bet,\bp)$, define the
{\it intersection dimension} $\Idim V$ to be the maximal number of a
priori given generic distinct points of $\Sig$ lying in $\bn_*\hat
C$ for an element $\{\bn:\hat C\to\Sig\}\in V$.

\subsection{Dimension count and generic elements}\label{secdim} We have a projection
\begin{equation}\pr:{\mathcal V}_\Sig(D,g,\alp,\bet,\bp)\to|D|,\quad\{\bn:\hat
C\to\Sig,\hat \bp\}\mapsto\bn_*\hat C\ ,\label{eproj}\end{equation}
and can define $V_\Sig(D,g,\alp,\bet,\bp)\subset|D|$ to be the
$\pr$-image of the union of the components of ${\mathcal
V}_\Sig(D,g,\alp,\bet,\bp)$ of the (maximal) intersection dimension
$R_\Sig(D,g,\bet)$.

\begin{proposition}\label{pn1} In the notations of section \ref{sec-gs},
let $D\in\Pic(\Sig,E)$, an integer $g\ge0$, and vectors
$\alp,\bet\in\Z^\infty_+$ satisfy (\ref{en1}), and let ${\mathcal
V}_\Sig(D,g,\alp,\bet,\bp)\ne\emptyset$. Then
\begin{enumerate}\item[(1)] $R_\Sig(D,g,\bet)\ge0$, and each component $V$ of
${\mathcal V}_\Sig(D,g,\alp,\bet,\bp)$ satisfies
\begin{equation}\Idim V\le R_\Sig(D,g,\bet)\ ;\label{enov1}\end{equation}
furthermore, the elements $\{\bn:\hat C\to\Sig,\hat p\}\in{\mathcal
V}_\Sig(D,g,\alp,\bet,\bp)$ with smooth $\hat C$ and $\#(\bn(\hat
C)\cap E\setminus\bp)\le n\le\|\bet\|$ form a stratum ${\mathcal
V}_{\Sig,n}(D,g,\alp,\bet,\bp)
\subset{\mathcal
V}_\Sig(D,g,\alp,\bet,\bp)$ of intersection dimension
\begin{equation}\Idim{\mathcal V}_{\Sig,n}(D,g,\alp,\bet,\bp)\le-(K_\Sig+E)D+g-1+n\ ,\label{enov11}\end{equation}
\item[(2)] if $D\ne sD_0$ for any $s\ge2$ and
divisor class $D_0$ such that $-(K_\Sig+E)D_0=0$, and (\ref{enov1})
becomes an equality, then
\begin{enumerate}\item[(2i)] a generic element $\{\bn:\hat
C\to\Sig,\hat\bp\}\in V$ is an immersion, birational onto
$C=\bn(\hat C)$, where $C$ is a nodal curve, nonsingular along $E$;
\item[(2ii)]if $R_\Sig(D,g,\bet)>0$, then the family $V$ has no base
points outside $\bp$, and the generic curve $C=\bn(\hat C)$ crosses
any a priori given curve $C'\subset\Sig$ transversally outside
$\bp$; furthermore, $V_\Sig(D,g,\alp,\bet,\bp)$ is smooth at
$C$.\end{enumerate}\end{enumerate}
\end{proposition}

{\bf Proof}. We divide our argument into several parts: in Steps 1-3
we prove statement (1), in Steps 4-8 we prove statement (2).

\smallskip
{\it Step 1: Preliminaries}. Without loss of generality, assume that
${\mathcal V}_\Sig(D,g,\alp,\bet,\bp)$ is irreducible.

From now on and till Step 7 below, we suppose that, for a generic
element $\qquad\qquad$ $\{\bn:\hat C\to\Sig\}\in {\mathcal
V}_\Sig(D,g,\alp,\bet,\bp)$, the map $\bn:\hat C\to C=\bn(\hat C)$
is birational, {\it i.e.} this is the normalization map.

Then inequalities (\ref{enov1}) and (\ref{enov11}) are evident if
$D$ is a $(-1)$-curve as well as if $D=dL-d_1E_1-...-d_6E_6$ with
$d\le 2$ or with $d_i<0$ for some $i$. So, we suppose that
$D=dL-d_1E_1-...-d_6E_6$ with $d\ge 3$, $d_1,...,d_6\ge 0$.

Due to generic position of $\bp$, one has
$$
\Idim {\mathcal V}_\Sig(D,g,\alp,\bet,\bp)=\Idim {\mathcal
V}_\Sig(D,g,0,\alp+\bet,\emptyset)-\|\alp\|\ ,
$$ and, furthermore, for the properties listed in (2i), (2ii), a generic
element of ${\mathcal V}_\Sig(D,g,\alp,\bet,\bp)$ is generic for
${\mathcal V}_\Sig(D,g,0,\alp+\bet,\emptyset)$, too. Hence we can
let $\alp=0$ and $\bp=\emptyset$. To shorten notation, we write
(within the present proof) ${\mathcal V}_\Sig(D,g,\bet)$ for
${\mathcal V}_\Sig(D,g,0,\bet,\emptyset)$.

Inequality (\ref{enov1}) and its proof are completely analogous to
\cite[Propositions 2.1 and 2.2]{CH}: namely, the argument is
developed for curves on any algebraic surface, and its applicability
amounts to checking a number of sufficient numerical conditions,
verified in \cite{CH} for the case of the plane. In the following we
do not copy the reasoning of \cite{CH}, but go through all numerical
conditions and verify them in our setting.

\smallskip
{\it Step 2: Proof of (\ref{enov1})}. The computation of $\Idim
{\mathcal V}_\Sig(D,g,\bet)$ literally goes along \cite[Section
2.3]{CH}, where one has to verify the following. \begin{enumerate}
\item[(D1)] The conclusion of \cite[Corollary 2.4]{CH} reads
in our situation as $$\Idim {\mathcal V}_\Sig(D,g,(DE)e_1)\le
-DK_\Sig+g-1\ ,$$ and it holds, since the hypothesis of
\cite[Corollary 2.4]{CH}, equivalent to $DK_\Sig<0$ is true for any
effective divisor class $D\in\Pic(\Sig)$.
\item[(D2)] The inequality
$\deg(\bn^*{\mathcal O}_{\PP^2}(1)(-\bd))\ge0$ 
in \cite[Page
363]{CH}, where $$\bd=\sum_{i\ge1,\ 1\le j\le\alp_i}i\cdot\hat
p_{ij}+ \sum_{i\ge1,\ 1\le j\le\bet_i}(i-1)\cdot\hat q_{ij}$$ (cf.
(\ref{ediv})), reads in our situation as $DE\ge I\bet-\|\bet\|$, and
it holds since $I\bet=DE$ by (\ref{en1}).
\item[(D3)] The
inequality \begin{equation}\deg(c_1({\mathcal N}_{\hat
C}(-\bd))\otimes\omega^{-1}_{\hat C})> 0\label{enn5}\end{equation}
in \cite[Page 363]{CH} (${\mathcal N}_{\hat C}$ is the normal sheaf
on $\hat C$ and $\omega_{\hat C}$ is the dualizing bundle) reads in
our setting as
$$\displaylines{
-DK_\Sig+ 2g-2 -\deg\bd + 2 - 2g \cr
=-D(K_\Sig+E)+\|\bet\|=(d-d_6)+\|\bet\| > 0, }
$$
and it holds true since $d-d_6\ge 1$ as $D$ is represented by a
reduced, irreducible curve with $d\ge 3$.\end{enumerate} After that,
as in the end of the proof of \cite[Proposition 2.1]{CH}, we derive
$$\Idim {\mathcal V}_\Sig(D,g,\bet)\le\deg(c_1({\mathcal
N}_{\hat C}(-\bd))-g+1=R_\Sig(D,g,\bet)\ ,
$$
which completes the proof of (\ref{enov1}).

\smallskip
{\it Step 3: Proof of (\ref{enov11}).} Again we can assume that
${\mathcal V}_{\Sig,n}(D,g,\alp,\bet,\bp)$ is irreducible and that
$\{\bn:\hat C\to\Sig,\hat\bp\}$ is generic in ${\mathcal
V}_{\Sig,n}(D,g,\alp,\bet,\bp)$. By our assumption, $\bn$ takes
$\hat C$ birationally onto its image, then (\ref{enov11}) can be
proven in the same way as the inequality in \cite[Corollary
2.7]{CH}, and, moreover, the proof literally coincides with the last
paragraph of the proof of \cite[Corollary 2.7]{CH}.

\smallskip
{\it Step 4: Immersion}. The sufficient conditions for $\bn:\hat
C\to\Sig$ to be an immersion as $\{\bn:\hat
C\to\Sig,\hat\bp\}\in{\mathcal V}_\Sig(D,g,\bet)$ is generic are
\begin{enumerate}\item[(I1)] $\deg(c_1({\mathcal N}_{\hat C}(-\bd))\otimes\omega^{-1}_{\hat
C})\ge 2$, equivalent to $d-d_6+\|\bet\|\ge2$, for the immersion
away from $\bn^{-1}(C\cap E)$ (see \cite[First paragraph of the
proof of Proposition 2.2]{CH}); \item[(I2)] $\deg(c_1({\mathcal
N}_{\hat C}(-\bd))\otimes\omega^{-1}_{\hat C})\ge 4$, equivalent to
$d-d_6+\|\bet\|\ge4$ for the immersion at $\bn^{-1}(C\cap E)$ (see
\cite[Second paragraph of the proof of Proposition 2.2]{CH}).
\end{enumerate}

Suppose that a generic $C=\bn(\hat C)$ has a local singular branch
at $z\in\Sig\setminus E$. Then by B\'ezout's bound for the
intersection of $C$ with a line $l\in|L-E_6|$ passing through $z$,
we get $d\ge d_6+2$, a contradiction by (I1).

Suppose that a generic $C=\bn(\hat C)$ has a singular local branch
at $z\in E$. Then $\|\bet\|\ge1$, and by the above B\'ezout's bound
we get $d-d_6+\|\bet\|\ge3$. Thus, by (I2), it remains to analyze
the case $d-d_6+\|\bet\|=3$, that means: $d_6=d-2$, $\|\bet\|=1$,
and $z$ is a center of a unique local branch of multiplicity $2$,
{\it i.e.} a singularity $A_{2s}$, $s\ge1$. Since $\del(A_{2s})=s$,
by the genus formula we have
\begin{equation}s\le\frac{D^2+DK_\Sig+2}{2}-g=d-2-g-\sum_{i=1}^5\frac{d_i(d_i-1)}{2}\
.\label{egenus}\end{equation} We shall show that, in such a case,
the family of curves of genus $g$ in the linear system $|D|$, having
singularity $A_{2s}$ at $z$, and satisfying $(C\cdot
E)_z=2d-d_1-...-d_5$ ({\it i.e.} $\|\bet\|=1$), has dimension $\le
g$. This will be a contradiction to the generality of $C$, since,
allowing the point $z$ to move along $E$, we would get
$$\Idim{\mathcal V}_\Sig(D,g,\bet)\le g+1=R_\Sig(D,g,\bet)-1\ .$$

Suppose that \begin{equation}2d-d_1-...-d_5=(C\cdot
E)_z=2s+1\label{einter}\end{equation} (the maximal possible
intersection number). Let $\widetilde\pi:\widetilde\Sig\to\Sig$ be
the blow up of $z$ and $s$ more infinitely near points of $C$ at
$z$. Then the considered family of curves on $\Sigma$ goes to the
family of immersed curves of genus $g$ in the linear system $|C^*|$
(the upper asterisk denotes the strict transform) passing through
the point $z^*\in\widetilde\Sig\setminus E^*$. Since the
multiplicity of $C$ at $z$ and the next $(s-1)$ infinitely near
points is $2$, and at the last blown up point is $1$, we have by
(\ref{einter})
$$\deg(c_1(\widetilde{\mathcal
N}_{\hat
C}(-z^*)))-(2g-2)=-K_{\widetilde\Sig}C^*-1=3d-d_1-...-d_6-2s-2=1>0\
,$$ (here $\widetilde{\mathcal N}_{\hat C}$ is the normal bundle of
$\widetilde\bn:\hat C\to C^*\hookrightarrow\widetilde\Sig$) which
ensures (cf. \cite[Page 357]{CH}) that the considered family has
dimension
\begin{equation}\le\deg(c_1(\widetilde{\mathcal N}_{\hat C}(-z^*)))-g+1=g\ ,\label{edim1}\end{equation} and we are
done, provided (\ref{einter}) holds.

Suppose that \begin{equation}2d-d_1-...-d_5=(C\cdot E)_z=2k,\quad
1\le k\le s\ .\label{einter1}\end{equation} Combining this with
(\ref{egenus}), we immediately obtain that
$$\sum_{i=1}^5d_i(2-d_i)\ge4+2g+2(s-k)\ ,$$ and hence (up to renumbering of $d_1,...,d_5$)
$$d_1=...=d_4=1,\ 0\le d_5\le2,\quad g=0,\quad k=s\
.$$ Blowing up $z$ and $(s-1)$ further infinitely near points of
$(C,z)$, we transform the germ at $C$ of the considered family of
curves into the germ at $C^*$ of the family of rational, immersed
curves in $|C^*|$ which are quadratically tangent to the last blown
up exceptional divisor in a neighborhood of the point $z^*$.
Similarly we compute
$$\deg(c_1(\widetilde{\mathcal
N}_{\hat
C}(-z^*)))-(2g-2)=-K_{\widetilde\Sig}C^*-1=3d-d_1-...-d_6-2s-1=1>0\
,$$ and derive the bound (\ref{edim1}) for the dimension of the
considered family which completes the proof of the immersion
property.

\smallskip

{\it Step 5: Nodality and nonsingularity along $E$}. From now on we
suppose that  $\bn:\hat C\to\Sig$ is an immersion. We may also
assume that $C\cap E\ne\emptyset$, {\it i.e.}, $\|\bet\|>0$. Indeed,
otherwise, we either replace $E$ by another $(-1)$-curve, or, if $C$
does not meet any $(-1)$-curve, we blow them down and reduce the
problem to the planar case.

In what follows, we argue by contradiction. Namely, assuming that
one of the statements fails, we derive that necessarily $\Idim
{\mathcal V}_\Sig(D,g,\bet)<R_\Sig(D,g,\bet)$.

Suppose that $\bn$ takes $s\ge2$ points of the divisor $\bd$ to the
same point $z\in E$. Fixing the position of $z$ in $E$, we obtain a
subvariety $U\subset {\mathcal V}_\Sig(D,g,\bet)$ of dimension
$\Idim U\ge \Idim {\mathcal V}_\Sig(D,g,\bet)-1$. On the other hand,
the same argument as in the proof of \cite[Proposition 2.1]{CH}
gives $$\Idim U\le h^0(\hat C,{\mathcal N}_{\hat C}(-\bd-\bd')),$$
where $\bd'=\sum_{i,j\ge1}\hat q_{ij}$. Assuming that $c_1({\mathcal
N}_{\hat C}(-\bd-\bd'))\otimes\omega^{-1}_{\hat C}$ is positive on
$\hat C$ and applying \cite[Observation 2.5]{CH}, we get (cf.
\cite[Page 364]{CH}) $$\Idim {\mathcal V}_\Sig(D,g,\bet)\le\Idim
U+1\le h^0(\hat C,{\mathcal N}_{\hat C}(-\bd-\bd'))+1$$
$$=\deg(c_1({\mathcal
N}_{\hat C}(-\bd-\bd'))-g+2=-DK_\Sig+2g-2-\deg(\bd+\bd')-g+2$$
$$=-D(K_\Sig+E)+g+\|\bet\|-s<R_\Sig(D,g,\bet)\ .$$ The
positivity required above is equivalent to
\begin{equation}-DK_\Sig-\deg(\bd+\bd')=d-d_6
+\|\bet\|-s > 0\ ,\label{ejan6_4}\end{equation} that holds due to
$\|\bet\|\ge s$, $d>d_6$.

Suppose that $|\bn^{-1}(z)|=s\ge3$ for some point $z\in\Sig\setminus
E$. Fixing the position of this point, we obtain a subvariety
$U\subset {\mathcal V}_\Sig(D,g,\bet)$ of dimension $\Idim U\ge\Idim
{\mathcal V}_\Sig(D,g,\bet)-2$. On the other hand, by the same
arguments as above we have again the inequality (\ref{ejan6_4}) and
an upper bound $$\Idim U\le h^0(\hat C,{\mathcal N}_{\hat
C}(-\bd-\bd'))\ ,$$ where $\bd'=\bn^{-1}(z)$. In view of
(\ref{ejan6_4}), the bundle $c_1({\mathcal N}_{\hat
C}(-\bd-\bd'))\otimes\omega^{-1}_{\hat C}$ is positive on $\hat C$,
and thus applying \cite[Observation 2.5]{CH} we get
$$\Idim
{\mathcal V}_\Sig(D,g,\bet)\le\Idim U+2\le h^0(\hat C,{\mathcal
N}_{\hat C}(-\bd-\bd'))+2$$
$$=\deg(c_1({\mathcal
N}_{\hat C}(-\bd-\bd'))-g+3=-DK_\Sig+2g-2-\deg(\bd+\bd')-g+3$$
$$=-D(K_\Sig+E)+g+\|\bet\|-s+1<R_\Sig(D,g,\bet)\ .$$

Suppose that $\bn^{-1}(z)=w_1+w_2$, $w_1\ne w_2\in\hat C$ for some
point $z\in\Sig\setminus E$, and the two local branches of
$C=\bn(\hat C)$ at $z$ intersect with multiplicity $s\ge2$. In
suitable coordinates in a neighborhood of $z$, $C$ is given by an
equation $y^2+2yx^s=0$, and the tangent space to the local
equisingular stratum in ${\mathcal O}_{\Sig,z}$ is the ideal
$I=\langle y+x^s,x^{s-1}y\rangle$, and it does not contain $y$. On
the other hand, the inequality $d-d_6+\|\bet\|>2$, which comes from
$\|\bet\|>0$ and the B\'ezout's bound $d\ge d_6+2$ for the
intersection of $C$ with the line $l\in|L-E_6|$ through $z$, can be
rewritten as $\deg({\mathcal N}_{\hat C}(-\bd-w_1-w_2))>2g-2$, and
hence $H^1(\hat C,{\mathcal N}_{\hat C}(-\bd-w_1-w_2))=0$. The
latter relation immediately implies that there is a section of the
bundle ${\mathcal N}_{\hat C}(-\bd)$ which vanishes at $w_1$ and
does not vanish at $w_2$. Moving $\bn:\hat C\to\Sig$ inside
${\mathcal V}_\Sig(D,g,\bet)$ along this section, we realize a
(one-parameter) deformation of the local equation $y^2+2yx^s=0$ of
$C$, tangent to the element $y\in{\mathcal O}_{\Sig,z}$, and hence
breaking the equisingularity. This contradiction proves that a
generic $C=\bn(\hat C)$ cannot have tangent local branches.

\smallskip
{\it Step 6: Absence of base points}. Assuming that
$R_\Sig(D,g,\bet)>0$ and ${\mathcal V}_\Sig(D,g,\alp,\bet,\bp)$ has
a base point $z\in\Sig\setminus\bp$, we obtain $$\Idim {\mathcal
V}_\Sig(D,g,\alp,\bet,\bp)\le h^0(\hat C,{\mathcal N}_{\hat
C}(-\bd-w))\ ,$$ where $\bn(w)=z$. To arrive to a contradiction, we
use the inequality
$$\Idim {\mathcal V}_\Sig(D,g,\alp,\bet,\bp)\le h^0(\hat C,{\mathcal
N}_{\hat C}(-\bd-\bd'))$$
$$=-DK_\Sig+2g-2-\deg(\bd+\bd')-g+1$$
$$=-D(K_\Sig+E)+g+\|\bet\|-2<R_\Sig(D,g,\bet)\ ,
$$ For it we need the positivity of the bundle $c_1({\mathcal
N}_{\hat C}(-\bd-w))\otimes\omega^{-1}_{\hat C}$ on $\hat C$, which
reduces to the inequality $d-d_6+\|\bet\|>1$. The only case it does
not hold is $d_6=d-1$, $\|\bet\|=0$, but then $g=0$ and
$R_\Sig(D,0,0)=d-d_6-1=0$, contrary to the assumptions made.
Furthermore, by Bertini we get the transversality of intersection of
$C=\bn(\hat C)$ with any a priori given curve $C'$ (outside $\bp$).

\smallskip
{\it Step 7: Smoothness}. The smoothness of
$V_\Sig(D,g,\alp,\bet,\bp)$ at $C$ is equivalent to
$$H^1(\hat C,{\mathcal N}_{\hat C}(-\bd))=0\quad\Longleftrightarrow\quad d-d_6+\|\beta\|>0\ ,$$
which clearly holds in the considered situation.

\smallskip
{\it Step 8: Multiple covers}. Assume that $\bn:\hat C\to C=\bn(\hat
C)$ is an $s$-fold covering, $s\ge2$, for a generic element
$\{\bn:\hat C\to\Sig,\hat\bp\}\in{\mathcal
V}_\Sig(D,g,\alp,\bet,\bp)$. Then $D=sD_0$, $2-2g=s(2-2g')-r$, where
$g'$ is the geometric genus of $C$, $r$ is the total ramification
multiplicity. Using (\ref{enov11}) for the normalization map
$\bn':C^{\bn}\to C\hookrightarrow\Sig$, we get (cf. \cite[Page
366]{CH} and \cite[Page 62]{Va})
$$\Idim {\mathcal V}_{\Sig,n}(D,g,\alp,\bet,\bp)\le-(K_\Sig+E)D_0+g'-1+\#(C\cap
E\setminus\bp)$$
\begin{equation}\le-\frac{(K_\Sig+E)D}{s}+\frac{g-1-r/2}{s}+n\le-(K_\Sig+E)D+g-1+n\
,\label{enov111}\end{equation} the latter inequality coming from the
fact that $-(K_\Sig+E)D\ge0$ for all effective divisor classes $D$.
Thus, we have proven both (\ref{enov1}) and (\ref{enov11}).

Now, the equality $\Idim {\mathcal
V}_\Sig(D,g,\alp,\bet,\bp)=R_\Sig(D,g,\bet)$, means the equality in
(\ref{enov111}) with $n=\|\bet\|$, and this leaves the only
possibility $(K_\Sig+E)D_0=0$, equivalent to $d=d_6$, and hence
either $D_0=E_i$, $1\le i\le5$, or $D_0=L-E_i-E_6$, $1\le i\le5$, or
$D_0=L-E_6$. If $D_0=E_i$, $1\le i\le5$, or $D_0=L-E_i-E_6$, $1\le
i\le5$, then $R_\Sig(D,g,\bet)=0$. If $D_0=L-E_6$ and
$R_\Sig(D,g,\bet)>0$, then $C=\bn(\hat C)$ is a generic line in the
one-dimensional linear system $L-E_6$, which, of course, crosses
transversally any a priori given curve $C'$.

The proof of Proposition \ref{pn1} is completed. \proofend

\begin{proposition}\label{pnov2}
Let $D\in\Pic_+(\Sig,E)$, and let $V$ be a component of a nonempty
family ${\mathcal V}_\Sig(D,g,\alp,\bet,\bp)$ such that, for a
generic element $\{\bn:\hat C\to\Sig,\hat\bp\}\in V$, $\bn:\hat C\to
C=\bn(\hat C)$ is an $s$-multiple cover, $s\ge2$. Then
\begin{enumerate}\item[(i)] either $D=sE_i$, $1\le i\le 5$, $C=E_i$, $g=0$,
$\alp=0$, $\bp=\emptyset$, $\bet=e_s$, $R_\Sig(D,g,\bet)=0$, and
$\bn:\hat C\to E_i$ has at least two critical points, one of which
is the preimage of $E_i\cap E$ and has ramification index $s$;
\item[(ii)] or $D=s(L-E_6)$, $C\in|L-E_6|$ is one of the two lines
quadratically tangent to $E$, $g=0$, $\alp=0$, $\bp=\emptyset$,
$\bet=e_{2s}$, $R_\Sig(D,g,\bet)=0$, and $\bn:\hat C\to C$ has at
least two critical points, one of which is the preimage of $C\cap E$
and has ramification index $s$;
\item[(iii)] or $D=s(L-E_6)$, $C\in|L-E_6|$ is a line crossing $E$
transversally at two points, $g=0$, $\alp=e_s$, $\bp$ is one of the
points of $C\cap E$, $\bet=e_s$, $R_\Sig(g,\bet)=0$, and $\bn:\hat
C\to C$ has two critical points projected onto $C\cap E$;
\item[(iv)] or $D=s(L-E_6)$, $C\in|L-E_6|$ is a line crossing $E$
transversally at two points, $g=0$, $\alp=0$, $\bp=\emptyset$,
$\bet=2e_s$, $R_\Sig(D,g,\bet)=1$, and $\bn:\hat C\to C$ has two
critical points projected onto $C\cap E$.
\end{enumerate}
\end{proposition}

{\bf Proof}. The statement is straightforward from the consideration
in Step 8 of the proof of Proposition \ref{pn1}. \proofend

\begin{proposition}\label{ptn2}
Let $D\in \Pic_+(\Sig,E)$, $g\ge0$, and $\alp,\bet\in\Z^\infty_+$
satisfy (\ref{en1}), and let $\bp=(p_{ij})_{i\ge1,1\le j\le\alp_i}$
be a sequence of generic distinct points of $E$. Assume that
$R_\Sig(D,g,\bet)=0$, and that ${\mathcal
V}_\Sig(D,g,\alp,\bet,\bp)$ contains an element $\{\bn:\hat
C\to\Sig,\hat\bp\}$ with $\bn$ birational onto its image. Then
\begin{enumerate}\item[(i)] either $D=E_i$, where $1\le i\le5$, and $g=0$,
$\bet=e_1$, $\#{\mathcal V}_\Sig(D,0,0,\bet,\emptyset)=1$;
\item[(ii)] or $D=L-E_i-E_6$, where $1\le i\le 5$, and $g=0$, $\bet=e_1$,
$\#{\mathcal V}_\Sig(D,0,0,\bet,\emptyset)=1$; \item[(iii)] or
$D=L-E_6$ and $g=0$, $\bet=e_2$, $\#{\mathcal
V}_\Sig(D,0,0,\bet,\emptyset)=2$;
\item[(iv)] or $D=dL-d_1E_1-...-d_5E_5-(d-1)E_6$,
where $d\ge1$, $0\le d_1,...,d_5\le1$, $d_1+...+d_5<2d$, $g=0$,
$I\alp=2d-d_1-...-d_5$, and  $\#{\mathcal
V}_\Sig(D,0,\alp,0,\bp)=1$;
\item[(v)] for $\alp\ne0$, $\bet\ne0$, we have $D=L-E_6$, $g=0$,
$\alp=\bet=e_1$, and  $\#{\mathcal V}_\Sig(D,0,\alp,\bet,\bp)=1$.
\end{enumerate}
\end{proposition}

{\bf Proof}. Straightforward from the formula
$R_\Sig(D,g,\bet)=d-d_6+g+\|\bet\|-1=0$. \proofend

\subsection{Statement of the problem}\label{sec23}
\begin{definition}\label{d23} Given a divisor class $D\in\Pic(\Sig,E)$, $g\ge0$, and vectors
$\alp,\bet\in\Z^\infty_+$ satisfying (\ref{en1}), and a sequence
$\bp$ of $\|\alp\|$ generic points of $E$, we define
$$N_\Sig(D,g,\alp,\bet)=\begin{cases}0,\quad &
\text{if}\ V_\Sig(D,g,\alp,\bet,\bp)=\emptyset,\\ \deg
V_\Sig(D,g,\alp,\bet,\bp),\quad & \text{if}\
V_\Sig(D,g,\alp,\bet,\bp)\ne\emptyset, \end{cases}$$
($V_\Sig(D,g,\alp,\bet,\bp)$ introduced in section
\ref{secdim}).\end{definition}

This number, of course, does not depend on the choice of $\bp$ and
equals the number of the curves in $V_\Sig(D,g,\alp,\bet,\bp)$
matching a generic configuration of $R_\Sig(D,g,\bet)$ points in
$\Sig\setminus E$. In the case $\alp=0$, $\bp=\emptyset$, and
$\bet=(DE)e_1$ ({\it i.e.} no assigned conditions on $E$), the
number $N_\Sig(D,g,0,(DE)e_1)$ coincides with the genus $g$
Gromov-Witten invariant $GW_g(\PP^2_6,D)$. So, we state

\smallskip

{\bf Problem}. {\it Compute all the numbers $N_\Sig(D,g,\alp,\bet)$
and, in particular, the Gromov-Witten invariants
$GW_g(\PP^2_6,D)=N_{\PP^2_6}(D,g,0,(DE)e_1)$}.

\smallskip

The answer is given by a Caporaso-Harris type recursive formula
(\ref{eq:main}) in Theorem \ref{t1} and the following initial values
for that formula, which can be derived from Propositions \ref{pnov2}
and \ref{ptn2}:

\begin{proposition}\label{pini}
Let $D\in \Pic_+(\Sig,E)$, $g\ge0$, and $\alp,\bet\in\Z^\infty_+$
satisfy (\ref{en1}). Then
\begin{enumerate}
\item[(1)] $N_\Sig(s(L-E_6),0,0,2e_s)=1$ for all $s\ge 1$;
\item[(2)] If $R_\Sig(D,g,\bet)=0$, then $N_\Sig(D,g,\alp,\bet)=0$
except for the following cases:
\begin{enumerate}\item[(i)] $N_\Sig(sE_i,0,0,e_s)=1$ for all
$s\ge1$, $i=1,...,5$; \item[(ii)] $N_\Sig(s(L-E_6),0,0,e_{2s})=2$
for all $s\ge1$;
\item[(iii)] $N_\Sig(s(L-E_6),0,e_s,e_s)=1$ for all $s\ge1$;
\item[(iv)] $N_\Sig(s(L-E_i-E_6),0,0,e_s)=1$ for all $s\ge1$, $i=1,...,5$;
\item[(v)] $N_\Sig(dL-d_1E_1-...-d_5E_5-(d-1)E_6,0,\alp,0)=1$ for all
$d\ge1$, $0\le d_1,...,d_5\le1$, $d_1+...+d_5<2d$,
$I\alp=2d-d_1-...-d_5$.
\end{enumerate}\end{enumerate}
\end{proposition}

\section{Degeneration}

Let $D\in\Pic(\Sig,E)$, $g\ge0$, and $\alp,\bet\in\Z^\infty_+$
satisfy relations (\ref{en1}) and the inequality
$R_\Sig(D,g,\bet)>0$. Let $\bp=(p_{ij})_{i\ge1,1\le j\le\alp_i}$ be
a sequence of $\|\alp\|$ generic points of $E$, and let $p$ be a
generic point in $E\setminus\bp$. Assume that ${\mathcal
V}_\Sig(D,g,\alp,\bet,\bp)\ne\emptyset$ and introduce
$${\mathcal V}^p_\Sig(D,g,\alp,\bet,\bp)=\big\{\{\bn:\hat
C\to\Sig,\hat\bp\}\in{\mathcal V}_\Sig(D,g,\alp,\bet,\bp)\ :\
p\in\bn(\hat C)\big\}\ .$$

\begin{proposition}\label{pdeg1}
Let $V$ be a component of ${\mathcal V}_\Sig(D,g,\alp,\bet,\bp)$ of
dimension $R_\Sig(D,g,\bet)$, and let $W$ be a component of
$V\cap{\mathcal V}^p_\Sig(D,g,\alp,\bet,\bp)$ of dimension
$R_\Sig(D,g,\bet)-1$. Then a generic element $\{\bn:\hat
C\to\Sig,\hat\bp\}$ of $W$ is as follows.

(1) Assume that $\bn(\hat C)$ does not contain $E$. Then $\hat C$ is
smooth, and there is $k\ge 1$ such that $W$ is a component of
${\mathcal
V}_\Sig(D,g,\alp+e_k,\bet-e_k,\bp\cup\{p_{k,\alp_k+1}=p\})$.

(2) Assume that $\bn(\hat C)$ contains $E$. Then
\begin{enumerate}\item[(i)] $\hat C=\hat E\cup X\cup Y\cup Z$, where $\bn:\hat E\to E$
is an isomorphism, $X$ is the union of components mapped by $\bn$ to
curves positively intersecting with $E$ and not belonging to the
linear system $|L-E_6|$, $Y$ is the union of components mapped by
$\bn$ to lines belonging to the linear system $|L-E_6|$, $Z$ is the
union of components contracted by $\bn$ to points;
\item[(ii)] $X=\bigcup_{1\le i\le k}\hat C^{(i)}$, $Y=\bigcup_{k<i\le m}\hat C^{(i)}$,
where $0\le k\le m$, and $\hat C^{(1)},...,\hat C^{(m)}$ are
disjoint smooth curves; \item[(iii)] for each $i=1,...,m$, the map
$\bn:\hat C^{(i)}\to\Sig$ represents a generic element in a
component of some space ${\mathcal
V}_\Sig(D^{(i)},g^{(i)},\alp^{(i)},\bet^{(i)},\bp^{(i)})$ of
intersection dimension $R_\Sig(D^{(i)},g^{(i)},\bet^{(i)})$ and such
that

- $\sum_{i=1}^mD^{(i)}=D-E$,

-
$\sum_{i=1}^mR_\Sig(D^{(i)},g^{(i)},\bet^{(i)})=R_\Sig(D,g,\bet)-1$,

- $\bp^{(i)}$, $i=1,...,m$, are disjoint subsets of $\bp$,

- $\bn:X\to\Sig$ is birational onto its image,

- the (reduced) curve $\bn(X\cup Y)$ is nodal  and nonsingular along
$E$,

- any quadruple $(D',g',0,\bet')$ with $R_\Sig(D',g',\bet')=0$
appears in the list $$(D^{(i)},g^{(i)},\alp^{(i)},\bet^{(i)}),\quad
1\le i\le k\ ,$$ at most once; \item[(iv)] if, for some
$i=k+1,...,m$, $\bn(\hat C^{(i)})=L'\in|L-E_6|$, where $L'$ is a
line tangent to $E$, then $\bn:\hat C^{(i)}\to L'$ is an
isomorphism.
\end{enumerate}
\end{proposition}

\begin{remark}\label{rmult} As a consequence, we deduce that
the element $\{\bn:\hat C\to\Sig,\hat \bp\}$ as in Proposition
\ref{pdeg1}(2) can be uniquely restored from $C=\bn_*\hat C$.
\end{remark}

{\bf Proof of Proposition \ref{pdeg1}(1,2i,2ii,2iii)}. We follow the
lines of \cite[Section 3]{CH}, where a similar statement
(\cite[Theorem 1.2]{CH}) is proven for the planar case, and also use
some arguments from the proof of \cite[Theorem 5.1]{Va}. Again we do
not copy all the details, but explain the main steps and perform all
necessary computations, which, in fact, are consequences of the
bounds (\ref{enov1}) and (\ref{enov11}), and of the relation
$(K_\Sig+E)E=-2$. The truly new claim, which we provide with a
complete proof in the very end, is that, in the case (2), $\bn_*\hat
C$ does not contain multiple $(-1)$-curves.

\smallskip

Consider a generic one-parameter family in the component $V$ of
${\mathcal V}_\Sig(D,g,\alp,\bet,\bp)$ with the central fibre
belonging to $W$. Using the construction of \cite[Section 3]{CH},
one can replace the given family with a family $\{\bn_t:\hat
C_t\to\Sig,\hat\bp_t\}_{t\in(\C,0)}$ having the same generic fibres
and a semistable central fibre such that (cf. conditions (b), (c),
(e) in\cite[Section 3.1]{CH}\footnote{We do not use conditions (a)
and (d) in \cite[Section 3.1]{CH}, which, in fact, are not
needed.}):
\begin{itemize}\item the family is represented by a surface
${\mathcal C}$ with at most isolated singularities and two morphisms
$\pi_\Sig:{\mathcal C}\to\Sig$, $\pi_\C:{\mathcal C}\to(\C,0)$, so
that for each $t\ne 0$, $\pi_\Sig:{\mathcal C}_t\overset{\rm
def}{=}\pi_\C^{-1}(t)\to\Sig$ is isomorphic to $\bn_t:\hat
C_t\to\Sig$, \item the fibre $\hat C_0={\mathcal C}_0$ is a nodal
curve,
\item the family is minimal with respect to the above properties.
\end{itemize} Notice that we have disjoint sections \begin{equation}t\in(\C,0)\setminus\{0\}\mapsto\hat p_{ij,t},\ i\ge1,\ 1\le
j\le\alp_i,\quad t\in(\C,0)\setminus\{0\}\mapsto\hat q_{ij,t},\
i\ge1,\ 1\le j\le\bet_i\ ,\label{esec1}\end{equation} defined by
(\ref{ediv}) for each fibre: \begin{equation}\bn_t^*(E\cap\bn_t(\hat
C_t))=\sum_{i\ge1,\ 1\le j\le\alp_i}i\cdot\hat
p_{ij,t}+\sum_{i\ge1,\ 1\le j\le\bet_i}i\cdot\hat q_{ij,t}\
,\label{esec2}\end{equation}
$$\bn_t(\hat p_{ij,t})=p_{ij}\in\bp,\ i\ge1,\ 1\le j\le\alp_i\
,$$ which close up at $t=0$ into some global sections, from which
$\hat p_{ij,0}$, $i,j\ge1$, remain disjoint.

\smallskip

In the case (1), if $\hat C_0$ splits into components $\hat
C^{(1)},...,\hat C^{(m)}$, $m\ge 1$, mapped to curves and some
components mapped to points, then the maps $\bn_0:\hat
C^{(i)}\to\Sig$ represent elements in some ${\mathcal
V}_\Sig(D^{(i)},g^{(i)},\alp^{(i)},\bet^{(i)},\hat\bp^{(i)})$,
$i=1,...,m$, for which we have
$$\sum_{i=1}D^{(i)}=D,\quad\sum_{i=1}^m\|\bet^{(i)}\|\le\|\bet\|-1,\quad\sum_{i=1}^m(g^{(i)}-1)\le
g-m$$ (since some section $q_{kj,t}$ closes up at $p$, and at least
$m-1$ intersection points are smoothed out when deforming $\hat C_0$
into $\hat C_t$, $t\ne0$). Hence (cf. \cite[Page 66]{Va})
$$\sum_{i=1}^m\Idim{\mathcal
V}_\Sig(D^{(i)},g^{(i)},\alp^{(i)},\bet^{(i)},\hat\bp^{(i)})=\sum_{i=1}^m(-(K_\Sig+E)D^{(i)}+g^{(i)}+\|\bet^{(i)}\|-1)$$
$$\le-(K_\Sig+E)D+d+\|\bet\|-1-m=R_\Sig(D,g,\bet)-m\ ,$$ which
implies $m=1$ and thus, in view of Proposition \ref{pn1}, the
statement (1).

\smallskip

In the case (2), assume that $\hat C_0=\hat E\cup\hat
C^{(1)}\cup...\cup\hat C^{(m)}\cup Z$, where $\bn_0(\hat E)=E$, the
components $\hat C^{(i)}$, $1\le i\le m$, are mapped by $\bn_0$ to
curves, and $\bn_0(Z)$ is finite. We have:
\begin{itemize}\item $\hat E$ splits into components $\hat E_1,...,\hat E_a$ such that
$(\bn_0)_*\hat E_i=s_iE$, $s_i\ge1$, $i=1,...,a$;
\item $\{\bn_0:\hat C^{(i)}\to\Sig,\hat\bp^{(i)}\}\in{\mathcal
V}_\Sig(D^{(i)},g^{(i)},\alp^{(i)},\bet^{(i)},\bp^{(i)})$, where
$\bp^{(i)}=\bn_0(\hat\bp^{(i)})\subset\bp$, and $\hat\bp^{(i)}=(\hat
p_{kj,0}\in\hat C^{(i)})_{k,j\ge1}$; furthermore,
$\bet^{(i)}=\widetilde\gam^{(i)}+\gam^{(i)}$, where $\gam^{(i)}$
labels the multiplicities of the sections $\hat q_{kj,t}$ which land
up on $\hat C^{(i)}$ for $t=0$, $i=1,...,m$; \item denote also by
$g'$ the sum of the geometric genera of all the algebraic components
of $Z$, and by $n'$ the number of those connected components $Z'$ of
$Z$, for which $\#Z'\cap\Sing(\hat C_0)$ is greater than the number
of algebraic components of $Z'$ plus $\#\{\hat z\in
Z'\cap\bigcup_i\hat C^{(i)}\ :\ \bn_0(\hat z)\in E\}$.
\end{itemize} By construction, the local branches $\bn_0:(\hat C^{(i)},\hat p_{kj,0})\to\Sig$
and $\bn_0:(\hat C^{(i)},\hat q_{kj,0})\to\Sig$ are deformed
continuously into respective branches $\bn_t:(\hat C_t,\hat
p_{kj,t})\to\Sig$ and $\bn_t:(\hat C_t,\hat q_{kj,t})\to\Sig$,
$t\ne0$. In turn, the other local branches $\bn_0:(\hat C^{(i)},\hat
z)\to\Sig$ with $\bn_0(\hat z)=z\in E$ are deformed so that they
become disjoint from $E$ as $t\ne0$; hence such a point $\hat z\in
\hat C^{(i)}$ must be either an intersection point of $\hat C^{(i)}$
with $\hat E$, or an intersection point of $\hat C^{(i)}$ with a
connected component of $Z$ which joins $\hat C^{(i)}$ with $\hat E$,
and, in the deformation $(\bn_t:\hat C_t\to\Sig)_{t\in(\C,0)}$, this
intersection point is smoothed out. Moreover, all nodes of $\hat
C_0$ are smoothed out in the deformation to $\hat C_t$, $t\ne0$.
Thus, $$\sum_{i=1}^a(g(\hat
E_i)-1)+\sum_{i=1}^m(g^{(i)}-1)+g'+n+n'+\sum_{i=1}^m\|\widetilde\gam^{(i)}\|\le
g-1\ ,$$ where $n$ is the number of nodes of $\hat C^{(1)},...,\hat
C^{(m)}$. Applying (\ref{enov11}), we get $$\Idim
W=R_\Sig(D,g,\bet)-1=-(K_\Sig+E)D+g-1+\|\bet\|-1$$
$$\le\sum_{i=1}^m\left(-(K_\Sig+E)D^{(i)}+
g^{(i)}-1+\#(\bn_0(\hat C^{(i)})\cap E\setminus\bp)\right)$$
$$=-(K_\Sig+E)(D-\sum_{i=1}^as_iE)+\sum_{i=1}^m(g^{(i)}-1)+\sum_{i=1}^m\|\bet^{(i)}\|$$
\begin{equation}\le-(K_\Sig+E)D+g-1-\sum_{i=1}^a(2s_i-1+g(\hat
E_i))+\sum_{i=1}^m\|\widetilde\gam^{(i)}\|-n-n'-g'\
.\label{enov2}\end{equation} Thus, we immediately obtain that
\begin{itemize}\item $a=1$, $s_1=1$, and $g(\hat E_1)=0$, {\it i.e.}
$\bn_0:\hat E\to E$ is an isomorphism; \item
$\sum_{i=1}^m\widetilde\gam^{(i)}=\bet$, {\it i.e.} all the points
$\hat q_{ij,0}$ belong to $\hat C^{(1)}\cup...\cup\hat C^{(m)}$;
\item for each $i=1,...,m$, $\{\bn_0:\hat C^{(i)}\to\Sig,\hat
\bp^{(i)}\}$ is a generic element of a component of ${\mathcal
V}_\Sig(D^{(i)},g^{(i)},\alp^{(i)},\bet^{(i)},\bp^{(i)})$ of
intersection dimension $R_\Sig(D^{(i)},g^{(i)},\bet^{(i)})$; in
particular, all the sections $t\in(\C,0)\mapsto\hat p_{kj,t}$ and
$t\in(\C,0)\mapsto\hat q_{kj,t}$, $k,j\ge1$, are pairwise disjoint
(cf. condition (d) in \cite[Section 3.1]{CH});
\item $n=0$,
{\it i.e.} $\hat C^{(1)},...,\hat C^{(m)}$ are disjoint.
\item $\|\gam^{(i)}\|>0$ for each $i=1,...,m$; \item $g'=0$, {\it
i.e.} all the algebraic components of $Z$ are rational, and $n'=0$,
{\it i.e.} the connected components of $Z$ are chains of rational
curves, which either join the points $\{\hat z\in \bigcup_i\hat
C^{(i)}\ :\ \bn_0(\hat z)\in E,\ \hat z\ne\hat p_{kj,0},\hat
q_{kj,0},\ k,j\ge1\}$ with $\hat E$ (cf. Assumption (b) in
\cite[Section 3.2]{CH}), or are attached by precisely one point to
$\hat E$ or to $\bigcup_i\hat C^{(i)}$; however, the latter type
connected components should not exist due to the minimality of
${\mathcal C}$.
\end{itemize} In view of Propositions \ref{pn1}, \ref{pnov2}, and \ref{ptn2}, it remains to prove that the points $\hat p_{kj,0}$,
which do not belong to $\hat C^{(1)}\cup...\cup\hat C^{(m)}$, lie on
$\hat E$. Indeed, if some point $\hat p_{kj,0}$ were in
$Z\setminus(\hat E\cup\bigcup_i\hat C^{(i)})$, then some point $\hat
z\in\bigcup_i\hat C^{(i)}\setminus\hat\bp_0$ would have been mapped
to $p_{kj}$, which would have resulted in reducing $1$ in the third
and fourth line of bounds (\ref{enov2}), thus, a contradiction.

\smallskip

The final step is to confirm that $C=(\bn_0)_*\hat C_0$ does not
contain multiple $(-1)$-curves. Let, for instance, $E_1$ have
multiplicity $s\ge2$ in $C$, and (in the above notations) let $\hat
C^{(k+1)},...,\hat C^{(m)}$ be all the components of $\hat C_0$
mapped onto $E_1$. Since $R_\Sig(D,g,\bet)>0$, we have
$\quad\quad\quad\quad$ $D=dL-d_1E_1-...-d_6E_6$, $d\ge 1$,
$d_1,...,d_5\ge 0$. Thus, $C'=(\bn_0)_*(\hat C^{(1)}\cup...\cup \hat
C^{(k)})$ belongs to the linear system
$|(d-2)L-(d_1+s-1)E_1-(d_2-1)E_2-...-(d_5-1)E_5-d_6E_6|$. That is,
$C'$ crosses $E_1\setminus E$ with multiplicity $d_1+s-1$, and as
explained above these intersection points persist in the deformation
$(\bn_t:\hat C_t\to\Sig)_{t\in(\C,0)}$. Hence, $(\bn_t)_*\hat C_t$
must cross $E_1$ with multiplicity $\ge d_1+s-1>d_1$, thus, a
contradiction. \proofend

The proof of Proposition \ref{pdeg1}(2iv) will be given in the end
of section \ref{sec542}. In turn the statement of Proposition
\ref{pdeg1}(2iv) will not be used before that.

\section{Deformation}\label{sec5}

\subsection{Preliminaries}\label{sec1}
Having generic elements of ${\mathcal V}^p_\Sig(D,g,\alp,\bet,\bp)$
as described in Proposition \ref{pdeg1}, we intend to deform them
into generic elements of ${\mathcal V}_\Sig(D,g,\alp,\bet,\bp)$ and
compute the following multiplicities.

Let $D\in\Pic_+(\Sig,E)$, a non-negative integer $g$, and vectors
$\alp,\bet\in \Z_+^\infty$ satisfy (\ref{en1}). Pick a sequence $\bp
= \{p_{ij}\}_{i\geq 1, 1 \leq j \leq \alp_i}$ of $\|\alp\|$ generic
distinct points on $E$. Assume, in addition, that
\begin{equation}n:=R_\Sig(D,g,\bet)>0\quad\text{and}\quad D\ne s(L-E_6),\ s\ge1\ .\label{eassum}
\end{equation} Then, by Proposition \ref{pnov2}, the components of ${\mathcal V}_\Sig(D,g,\alp,\bet,\bp)$
of intersection dimension $R_\Sig(D,g,\bet)$ have the genuine
dimension $R_\Sig(D,g,\bet)$, and their union is birational to its
image $V_\Sig(D,g,\alp,\bet,\bp)$ in $|D|$ by the projection
(\ref{eproj}). Pick a set $\overline\bp$ of $n-1$ generic points of
$\Sig\setminus E$. Then
\begin{equation}V_\Sig(D,g,\alp,\bet,\bp,\overline\bp)\overset{\text{def}}{=}\left\{C\in
V_\Sig(D,g,\alp,\bet,\bp)\ :\
C\supset\overline\bp\right\}\label{estratum}\end{equation} is
one-dimensional (or empty).

Let $p\in E\setminus\bp$ be a generic point, and let $V$ be a
component of ${\mathcal V}_\Sig^p(D,g,\alp,\bet,\bp)$ of
(intersection) dimension $R_\Sig(D,g,\bet)-1$. Let $\{\bn:\hat
C\to\Sig,\hat\bp\}\in V$ be such that $p\in C=\bn_*\hat C$. Within
the given data, put $V(\overline\bp,C)=\emptyset$ if $C\not\in
V_\Sig(D,g,\alp,\bet,\bp,\overline\bp)$, and put $V(\overline\bp,C)$
to be the germ at $[C]$ of $V_\Sig(D,g,\alp,\bet,\bp,\overline\bp)$
if $[C]\in V_\Sig(D,g,\alp,\bet,\bp,\overline\bp)$. By construction,
$\{\bn:\hat C\to\Sig,\hat\bp\}$ is a generic element of $V$; hence
by Propositions \ref{pnov2} and \ref{pdeg1}, $C$ is nonsingular at
$p$. Take a smooth curve germ $L_p$ transversally crossing $E$ and
$C$ at $p$. Thus, we have a well defined map $\varphi_C:V(\overline
\bp,C)\to L_p$. The multiplicities we are interested in are the
degrees $\deg\varphi_C$ (equal $0$ if $V(\overline
\bp,C)=\emptyset)$.

\subsection{Deformation of an irreducible degenerate curve}

\begin{proposition}\label{p5}
In the notations and assumptions of section \ref{sec1}, let
$\bet_k>0$ and $\quad$ $\quad$ $\quad$ $\quad$ \mbox{$[C]\in
V_\Sig(D,g,\alp+e_k,\bet-e_k,\bp\cup\{p=p_{k,\alp_k+1}\})$}. Then
$\deg\varphi_C=k$.
\end{proposition}

{\bf Proof}. In view of (\ref{eassum}), $C$ possesses the properties
listed in Proposition \ref{pn1}(2). Furthermore, the germ of
$V_\Sig(D,g,\alp+e_k,\bet-e_k,\bp\cup\{p=p_{k,\alp_k+1}\})$ at $[C]$
is smooth by Proposition \ref{pn1}(2). Excluding $p$ from the set of
fixed points, we deduce that $V_\Sig(D,g,\alp,\bet,\bp)$ contains
$[C]$ and is smooth at $[C]$ as well. In particular,
$$h^0(\hat C,{\mathcal N}_{\hat C}(-\bd-p))=n-1\quad\text{and}\quad
h^0(\hat C,{\mathcal N}_{\hat C}(-\bd))=n\ .$$ In view of the
generic choice of $\overline\bp$, we obtain
$$h^0(\hat C,{\mathcal N}_{\hat C}(-\bd-p-\overline\bp))=0\quad\text{and}\quad
h^0(\hat C,{\mathcal N}_{\hat C}(-\bd-\overline\bp))=1\ ,$$ which in
its turn means that the germ $V(\overline\bp,C)$ intersect
transversally at $[C]$ with the hyperplane $H_p=\{[C']\in|D|\ :\
p\in C'\}$ in the linear system $|D|$. The latter conclusion yields
that $V(\overline\bp,C)$ diffeomorphically projects onto the germ of
$E$ at $p$ by sending a curve $[C']\in V(\overline\bp,C)$ to its
intersection point with $E$ in a neighborhood of $p$. Thus, in
suitable local coordinates $x,y$ of $\Sig$ in a neighborhood of $p$,
we have $p=(0,0)$, $E=\{y=0\}$, $L_p=\{x=0\}$, and a local
parametrization
\begin{equation}C_t=\left\{ay+b(x+t)^k+\sum_{i+kj\ge k}O(t)\cdot
(x+t)^iy^j=0\right\},\quad a,b\in\C^*,\ t\in(\C,0)\
,\label{e30}\end{equation} of $V(\overline\bp,C)$. Hence, for any
point $p'=(0,\tau)\in L_p$, we obtain precisely $k$ curves $C_t\in
V(\overline\bp,C)$ passing through $p'$ and corresponding to the $k$
values of $t$
\begin{equation}t=\left(-\frac{a}{b}\right)^{1/k}\tau^{1/k}+\
\text{h.o.t.}\ ,\label{e15}\end{equation} which completes the proof.
\proofend

\subsection{Deformation of a reducible degenerate curve}
We start with an auxiliary statement.

\begin{lemma}\label{l1} (1)
In the notations and assumptions of section \ref{sec1}, let
$C\supset E$. Then $C$ splits into distinct irreducible components
as follows \footnote{The coefficients $s',s''$ in the formula denote
the multiplicities of $L',L''$ in $C$, respectively.}:
\begin{equation}C=E\cup\bigcup_{i=1}^mC^{(i)}
\cup s'L'\cup s''L''\ ,\label{eC}\end{equation} where
\begin{enumerate}\item[(i)] $L',L''$ are the two lines in $|L-E_6|$
tangent to $E$, \item[(ii)] $C^{(1)},...,C^{(m)}$ do not contain
neither $L'$, nor $L''$, and, furthermore, each $C^{(i)}$ either is
a reduced, irreducible curve, or is $kL(p_{kl})$, where $k\ge 2$,
$p_{kl}\in\bp$, $L(p_{kl})\in|L-E_6|$ contains $p_{kl}$, or is
$kL(z)$, where $k\ge 2$, $z\in\overline\bp$, $L(z)\in|L-E_6|$
contains $z$,
\item[(iii)] the curve $C_{\redu}$ is nodal and nonsingular along
$E$.\end{enumerate} In addition, \begin{enumerate}\item[(iv)] for
each $i=1,...,m$, $C^{(i)}$ is a generic element in some
$V_\Sig(D^{(i)},g^{(i)},\alp^{(i)},\bet^{(i)},\bp^{(i)})$,
\item[(v)] $\sum_{i=1}^mD^{(i)}=D-E-(s'+s'')(L-E_6)$, \item[(vi)]
$\sum_{i=1}^mR_\Sig(D^{(i)},g^{(i)},\bet^{(i)})=R_\Sig(D,g,\bet)-1$,
\item[(vi)] $\bp^{(i)}$, $i=1,...,m$, are disjoint subsets of $\bp$;
\item[(vii)]
$\overline\bp^{(i)}:=\overline\bp\cap C^{(i)}$, $i=1,...,m$. form a
partition of $\overline\bp$;
\item[(viii)] there is a sequence of vectors $\gam^{(i)}\in\Z_+^\infty$, $i=1,...,m$ such that
\begin{equation}0<\gam^{(i)}\le\bet^{(i)}\quad\text{and}\quad
\bet=\sum_{i=1}^m(\bet^{(i)}-\gam^{(i)})\ .\label{eB}\end{equation}
\end{enumerate}

(2) Let $C$ be as above, and let $V(\overline\bp,C)\ne\emptyset$.
Let $\{C_t\}_{t\in(\C,0)}$ be a parameterized branch of
$V(\overline\bp,C)$ centered at $C=C_0$, and let $\{\bn_t:\hat
C_t\to\Sig,\hat\bp_t\}$, $t\in(\C,0)$ be its lift to ${\mathcal
V}_\Sig(D,g,\alp,\bet,\bp)$. Then in the deformation $C_0\to C_t$,
$t\ne0$, we have
\begin{enumerate}\item[(i)] a point $z\in\Sing(C_{\redu})\setminus E$ is either a node of some $C^{(i)}$, $1\le i\le m$, and then
it deforms into a moving node of $C_t$, $t\ne 0$, or $z$ is a
transverse intersection point of distinct irreducible components
$C',C''$ of $C_{\redu}$, both smooth at $z$, and then the point $z$
deforms into $k'k''$ nodal points of $C_t$, $t\ne0$, where $k'$ and
$k''$ are the multiplicities of $C'$, $C''$ in $C$, respectively;
\item[(ii)] if $z=p_{kl}\in\bp\setminus\bigcup_{i=1}^mC^{(i)}$,
then the germ of $E$ at $z$ turns into a smooth branch of $C_t$
crossing $E$ at $z$ with multiplicity $k$; \item[(iii)] if
$z=p_{kl}\in\bp\cap C^{(i)}$ for some $i=1,...,m$, then the germ of
$C^{(i)}$ at $z$ deforms into a smooth branch of $C_t$ centered at
$z$ and crossing $E$ with multiplicity $k$, and the germ of $E$ at
$z$ deforms into a smooth branch of $C_t$, $t\ne0$, disjoint from$E$
and crossing the former branch of $C_t$ transversally at $k$ points;
\item[(iv)] if $z=q_{kl}\in C^{(i)}\cap E\setminus\bp$ for some
$i=1,...,m$, and $q_{kl}=\bn(\hat q_{kl,0})$, where $(\hat
q_{kl,t})_{t\in(\C,0)}$ is defined by (\ref{esec1}), (\ref{esec2}),
then, first, in case of $C^{(i)}=sL(z)$, $z\in\overline\bp$, we have
$s=k$, and, second, the germ of $C^{(i)}$ at $z$ deforms into a
smooth branch of $C_t$, $t\ne0$, crossing $E$ at one point in a
neighborhood of $z$ with multiplicity $k$, and the germ of $E$ at
$z$ deforms into a smooth branch of $C_t$, $t\ne0$, disjoint from$E$
and crossing the former branch of $C_t$ transversally at $k$ points;
\item[(v)] if a point $z\in C^{(i)}\cap E\setminus\bp$ is not
$\bn(\hat q_{kl,0})$ for any section $(\hat q_{kl,t})_{t\in(\C,0)}$
as in (iv), then the germ of $E\cup C^{(i)}$ at $z$ deforms into an
immersed cylinder disjoint from $E$;
\item[(vi)] if $s'>0$ (or $s''>0$) then the union of $s'L'$
(resp., $s''L''$) with the germ of $E$ at the point $z=E\cap L'$
(resp., $z=E\cap L''$) turns into an immersed disc disjoint from $E$
and having $s'$ (resp., $s''$) nodes.
\end{enumerate}
\end{lemma}

{\bf Proof}. All the statements follow from Proposition
\ref{pdeg1}(2). \proofend

{\bf Notation.} In the sequel, the symbols
$C,C^{(i)},L(p_{kl}),L(z),L',L''$ always mean curves in $\Sig$ as
described in Lemma \ref{l1}.

\smallskip

Our main deformation result is:

\begin{proposition}\label{p6} Let us be given $D$, $g$, $\alp$, $\bet$,
$\bp$, $\overline\bp$ satisfying the assumptions of section
\ref{sec1}. Let $C\in|D|$ be given by (\ref{eC}) and satisfy the
conditions of Lemma \ref{l1}(1). Then
\begin{equation}\deg\varphi_{C}=(s'+1)(s''+1)\sum_{\{\gam^{(i)}\}_{i=1}^{m}}\prod_{i\in S}I^{\gam^{(i)}}\prod_{i=1}^m
\binom{\bet^{(i)}}{\gam^{(i)}}\ ,\label{e18}\end{equation} where the
sum runs over all sequences $\gam^{(i)}\in\Z_+^\infty$, $i=1,...,m$,
subject to (\ref{eB}), and $S$ is the set of those $i=1,...,m$ for
which $C^{(i)}$ is not of the form $kL(p_{kl})$.
\end{proposition}

\subsection{Proof of Proposition \ref{p6}}\label{sec2} Our strategy is as follows. Let
$$(L_p,p)\to(\C,0),\quad p'\in L_p\mapsto\tau\in(C,0)\ ,$$ be a parametrization of $L_p$.
Assume that $V(\overline\bp,C)\ne\emptyset$, take one of its
irreducible branches $V=\{C_t\ :\ t\in(\C,0)\}$ with a
parametrization normalized by a relation \begin{equation}p'=L_p\cap
C_t\quad\Longleftrightarrow\quad\tau=t^\mu\label{e22}\end{equation}
for some $\mu\ge1$. Clearly, $\deg\varphi_{C}\big|_{V}=\mu$, and the
curves $C\in V$ passing through $p'$ can be associated with
different parameterizations of $V$ obtained by substitutions
$t\mapsto t\eps$, $\eps^\mu=1$. Next, to a parameterized branch $V$
we assign a collection of bivariate polynomials (called
\emph{deformation patterns}) which describe deformations of the
curve $C$ in a neighborhood of isolated singularities of
$\bigcup_{i=1}^mC^{(i)}\cap E$ and in a neighborhood of non-reduced
components $C^{(i)}$, and in neighborhood of $L',L''$, when moving
along the branch $V$. We show that the constructed deformation
patterns belong to an explicitly described finite set, and, finally,
prove that, given a collection of arbitrary deformation patters from
the above sets, there exists a unique parameterized branch $V$ of
$V(\overline\bp,C)$ which induces the given deformation patterns.
Thus, the degree $\deg\varphi_{C}$ appears to be the number of
admissible collections of deformation patterns in the construction
presented below.

We have chosen deformation patterns in view of their convenience for
further real applications (see \cite{IKS1}), particularly, they suit
well for the calculation of Welschinger signs of real nodal curves.
Another advantage of deformation patterns is that we will be able to
complete the classification of degenerations in Proposition
\ref{pdeg1}, claim (2iv), by showing that the lines $L'$, $L''$
cannot be multiply covered.

In sections \ref{sec541} and \ref{sec542} we present a detailed
construction of deformation patterns for isolated singularities and
for multiple lines $L',L''$, which are the most involved cases. In a
similar way, in sections \ref{sec444} and \ref{sec445} we construct
deformation patterns for multiply covered lines $L^{\bp}_j$ and
$L_j$, but omit routine details.

\subsubsection{Preliminary choice in the construction of deformation
patterns}\label{sec441} First, we fix a sequence
$\gam^{(i)}\in\Z_+^\infty$, $i=1,...,m$, satisfying (\ref{eB}).

Next, for each $i=1,...,m$, take a subset of $\|\gam^{(i)}\|$ points
of $C^{(i)}\cap E\setminus\bp$ selected in such a way that, for
precisely $\gam^{(i)}_k$ of them, the intersection multiplicity of
$C^{(i)}$ and $E$ at the chosen point equals $k$. We notice that the
number of possible choices is
$\prod_{i=1}^m \left(\begin{matrix}\bet^{(i)}\\
\gam^{(i)}\end{matrix}\right)$, and we denote the chosen points by
$q'_{kj}$, $k\ge1$, $1\le j\le\sum_i\gam^{(i)}_k$, assuming that
$(C^{(i)}\cdot E)(q'_{kj})=k$ as $q'_{kj}\in C^{(i)}$. In
particular, for a component $C^{(i)}=kL(p_{kl})$, the only point of
$L(p_{kl})\cap E\setminus\bp$ is always chosen. Denote by $\bz'$ the
set of all the selected points.

In the next sections \ref{sec541}, \ref{sec542}, \ref{sec444}, and
\ref{sec445}, we suppose that $V(\overline\bp,C)\ne\emptyset$, and
there exists a branch $V$ of $V(\overline\bp,C)$ parameterized with
normalization (\ref{e22}) and which lifts to a family $\{\bn_t:\hat
C_t\to\Sig,\hat\bp_t\}\in{\mathcal V}_\Sig(D,g,\alp,\bet,\bp)$,
$t\in(\C,0)$, such that $(\bn_0)_*\hat C_0=C$, and the points
$\bn_0(\hat q_{ij,0})$ do not belong to $\bz'$ (here the sections
$\hat q_{ij,t}\in\hat C_t$ are uniquely defined by (\ref{esec1}),
(\ref{esec2})). Then put $\bz=\{\bn_0(\hat q_{ij,0})\ :\ i,j\ge1\}$.

\subsubsection{Deformation patterns for isolated
singularities}\label{sec541} Pick a point $q'_{kj}\in\bz'\cap
\bigcup_{i=1}^mC^{(i)}$.

Let $\Pi_E:\Sig\to\PP^2$ be the blow down of the (disjoint)
$(-1)$-curves
$$L-E_1-E_2,\ L-E_2-E_3,\ L-E_1-E_3,\ E_4,\ E_5,\ E_6\ .$$ In particular, it takes $E$
to the line through the points $\Pi_E(E_4)$ and $\Pi_E(E_5)$. Take
an affine plane $\C^2\subset\PP^2$ with the coordinates $u,v$ such
that $\Pi_E(E_4)=(0,0)$, $\Pi_E(E)\cap\C^2=\{v=0\}$. For the sake of
notation, we write $u(z),v(z)$ for the coordinates $u,v$ of the
point $\Pi(z)$ as $z$ is a point or a contracted divisor in $\Sig$.
For example, we write $p=(u(p),0)$ with some $u(p)\ne 0$, and,
furthermore, we suppose that $L_p=\{u=u(p)\}$.

Notice that this identification naturally embeds $H^0(\Sig,{\mathcal
O}_\Sig(D))$ into the space $P_{d_0}$ of polynomials in $u,v$ of
degree $\le d_0=2d-d_1-d_2-d_3$ (where $D=dL-d_1E_1-...-d_6E_6$).

The curve $C$ is then given by
\begin{equation}F_C(u,v):=v\left(f_1(u)+v\sum_{k,l\ge 0}a_{kl}u^kv^l\right)=0\ ,\label{e20}\end{equation} where
\begin{equation}f_1(u)=u^{d_4-1}\prod_{q'_{ij}\in\bz'}(u-u(q'_{ij}))^i\cdot\prod_{z\in\widetilde
C\cap E\setminus\bz'}(u-u(z))^{(\widetilde C'\cdot E)(z)}\
,\label{enn10}\end{equation} where $\widetilde C$ is the union of
the components of $C$ different from $E$.

\smallskip

A curve $C_t\in V$ is given by an equation
\begin{equation}at^\lambda
(f_2(u)+O(t))+v\left(f_1(u)+O(t)+v\sum_{k,l\ge
0}(a_{kl}+O(t))u^kv^l\right)=0\ ,\label{e24}\end{equation} where
\begin{equation}f_2(u)=u^{d_4}(u-u(E_5))^{d_5}\prod_{i\ge
1}\prod_{j=1}^{\alp_i}(u-u(p_{ij}))^i\cdot\prod_{i\ge
1}\prod_{j=1}^{\bet_i}(u-u(q_{ij}))^i\label{e21}\end{equation} (here
$p_{ij}$ runs over the points of $\bp$, and $q_{ij}$ are the
projections of the central points $\hat q_{ij,0}$ of the sections
$\hat q_{ij,t}$ defined by (\ref{esec1})). From (\ref{e22}), we
immediately derive that
\begin{equation}\lambda=\mu,\quad a=-\frac{f_1(u(p))}{f_2(u(p))}\
.\label{e23}\end{equation}

\smallskip

In the coordinates $x=u-u(q'_{ij})$, $y=v$, where $q'_{ij}\in\bz'$
is the chosen point, the curve $C$ is given by an equation
$$a_{02}^{ij}y^2+a_{i1}^{ij}x^iy+\sum_{k+il>2i}a^{ij}_{kl}x^ky^l=0,\quad a^{ij}_{02},a^{ij}_{i1}\ne 0\ ,$$
and the curve $C_t\in V$ ($t\ne 0$) is given by an equation
$$\Phi(x,y,t):=a^{ij}_{02}y^2+a^{ij}_{i1}x^iy+\sum_{k+il=2i}O(t)\cdot x^ky^l+\sum_{k+il>2i}(a^{ij}_{kl}+O(t))x^ky^l$$ $$+
\sum_{k+il<2i}t^{\bn(k,l)}(a^{ij}_{kl}+O(t))x^ky^l=0\ ,$$ where the
last sum does contain nonzero monomials, and for all of them
$\bn(k,l)>0$ and $a^{ij}_{kl}\ne0$, in particular, by (\ref{e24})
and (\ref{e23}),
\begin{equation}\bn(0,0)=\mu,\quad
a^{ij}_{00}=-\frac{f_1(u(p))f_2(u(q'_{ij}))}{f_2(u(p))}\
.\label{e27}\end{equation} Making an additional coordinate change
$$x\mapsto x-yt^{\bn(i-1,1)}\left(\frac{a^{ij}_{i1}}{a^{ij}_{02}}+O(t)\right)(a^{ij}_{i-1,1}+O(t))\ ,$$ we
can annihilate the monomial $x^{i-1}y$ in $\Phi(x,y,t)$.

Notice that
\begin{equation}\bn(k,0)\ge\bn(0,0)=\mu, \quad k>0\ .\label{e25}\end{equation} Indeed, otherwise,
we would have an intersection point of $C_t$ and $E$ converging to
$q'_{ij}$, a contradiction. Hence there exists
\begin{equation}\rho=\min_{k+il<2i}\frac{\bn(k,l)}{2i-k-il}>0\ .\label{e26}\end{equation} Now we
consider the equation
$\Psi(x,y,t):=t^{-2i\rho}\Phi(xt^\rho,yt^{i\rho},t)=0$. We have
$$\Psi(x,y,t)=a^{ij}_{02}y^2+a^{ij}_{i1}x^iy+\sum_{k+il=2i}O(t)\cdot x^ky^l+\sum_{k+il>2i}(a^{ij}_{kl}+O(t))t^{\rho(k+il-2i)}x^ky^l$$
$$+\sum_{k+il\le2i}t^{\bn(k,l)-\rho(2i-k-il)}(a^{ij}_{kl}+O(t))x^ky^l$$
and the well-defined polynomial
$$\Psi(x,y,0)=a_{02}^{ij}y^2+a_{i1}^{ij}x^iy+\sum_{k+il<2i}b^{ij}_{kl}x^ky^l$$ with the
coefficients $$b^{ij}_{kl}=\begin{cases}a^{ij}_{kl},\quad &
\text{if}\ \bn(k,l)=\rho(2i-k-il),\\ 0,\quad & \text{if}\
\bn(k,l)>\rho(2i-k-il)\ .\end{cases}$$ In view of (\ref{e25}),
$b_{k,0}=0$, $k\ge 1$, and hence
\begin{equation}\Psi(x,y,0)=a^{ij}_{02}y^2+a^{ij}_{i1}x^iy+\sum_{k=0}^{i-2}b^{ij}_{k1}x^ky+b^{ij}_{00}\
,\label{e39}\end{equation} where at least one of the $b^{ij}_{k1}$
or $b^{ij}_{00}$ is nonzero. We claim that $b^{ij}_{00}\ne0$.
Indeed, otherwise, $\Psi(x,y,t)=0$ would split into the line
$E=\{y=0\}$ and a curve, crossing $E$ in at least two points, which
in turn would mean that, in the deformation of $C$ into $C_t$ in a
neighborhood $U_{ij}$ of $q'_{ij}$, the germs of $\widetilde C$ and
$E$ at $q'_{ij}$ would glue up with the appearance of at least two
handles, thus breaking the claim of Lemma \ref{l1}(2v). Hence
$b^{ij}_{00}=a^{ij}_{00}$ given by (\ref{e27}), and as observed
above the equation $\Psi(x,y,0)=0$ defines an immersed affine
rational curve. We call $\Psi(x,y,0)$ the \emph{deformation pattern}
for the point $q'_{ij}\in\bz'$ associated with the parameterized
branch $V$ of $V(\overline\bp,C)$. Its Newton polygon is depicted in
Figure \ref{f1}. We summarize the information on deformation
patterns in the following statement:

\begin{figure}
\setlength{\unitlength}{1cm}
\begin{picture}(6,4)(-5,0.5)
\thicklines \put(1,1){\line(0,1){2}}\put(1,1){\line(3,1){3}}
\put(1,3){\line(3,-1){3}}
\thinlines\put(1,1){\vector(1,0){4}}\put(1,1){\vector(0,1){3}}\dashline{0.2}(1,2)(4,2)\dashline{0.2}(4,1)(4,2)
\put(3.9,0.5){$i$}\put(0.7,1.9){$1$}\put(0.7,2.9){$2$}
\end{picture}
\caption{Newton polygon of deformation pattern of an isolated
singularity}\label{f1}
\end{figure}
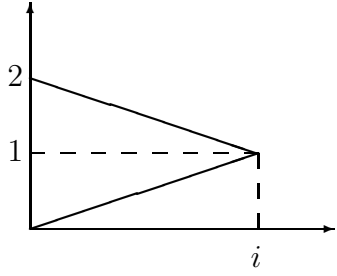

\begin{lemma}\label{l2}
(1) Given a curve $C=E\cup\widetilde C$ and a set $\bz'\subset
\widetilde C\cap E$ as chosen in section \ref{sec441}, the
coefficients $a^{ij}_{02}$, $a^{ij}_{i1}$, and
$b^{ij}_{00}=a^{ij}_{00}$ of a deformation pattern of the point
$q'_{ij}\in\bz'$ are determined uniquely up to a common factor.

(2) The set ${\mathcal P}(q'_{ij})$ of polynomials given by formula
(\ref{e39}), having fixed coefficients $a^{ij}_{02}$, $a^{ij}_{i1}$,
and $b^{ij}_{00}=a^{ij}_{00}$  and defining a rational curve,
consists of $i$ elements. Moreover, their coefficients
$b^{ij}_{k1}$, where $k\equiv i+1\mod2$, vanish, and, given a
polynomial $\Psi(x,y)\in{\mathcal P}(q'_{ij})$, the other members of
${\mathcal P}(q'_{ij})$ can be obtained by the following
transformations
\begin{eqnarray}&\Psi(x,y)\mapsto\Psi(x\eps,y),\ \text{where}\
\eps^i=1,\quad\text{if}\ i\equiv1\mod2\ ,\nonumber\\
&\Psi(x,y)\mapsto\Psi(x\eps,y\eps^i),\ \text{where}\
\eps^{2i}=1,\quad\text{if}\ i\equiv0\mod2\ .\nonumber\end{eqnarray}

(3) The germ at $F\in{\mathcal P}(q'_{ij})$ of the family of those
polynomials with Newton triangle $\conv\{(0,0),(0,2),(i,1)\}$, which
define a rational curve, is smooth of codimension $(i-1)$, and
intersect transversally (at one point $\{F\}$) with the affine space
of polynomials
$$\left\{a_{02}^{ij}y^2+a_{00}^{ij}+a_{i1}^{ij}x^iy+\sum_{k=0}^{i-2}c_kx^ky\
:\ c_0,...,c_{i-2}\in\C\right\}\ .$$
\end{lemma}

{\bf Proof}. The first statement comes from construction. In
particular, $a^{ij}_{00}$ is determined by formula (\ref{e27}). The
second statement follows from \cite[Lemma 3.5]{Sh0}. At last, the
second statement additionally implies that the family in the third
statement is the (germ of the) orbit of $F$ by the action
$$F(x,y)\mapsto \xi F(\xi_1x+\xi_0,\xi_2y),\quad\xi,\xi_1,\xi_2\in\C^*,\ \xi_0\in\C\ ,$$ and the third claim follows. \proofend

\begin{remark}\label{r5}
(1) Since the coefficients $b^{ij}_{k1}$ in (\ref{e39}) vanish for
$k\equiv i+1\mod2$, and the other ones not, the values
$$\bn(i-2k,1)=\rho(2i-(i-2k)-i)=\frac{\mu}{2i}\cdot 2k=\frac{\mu}{i},\quad 1\le k\le\frac{i}{2}\
,$$ are integer; hence $\mu$ is divisible by $i$.

(2) The treatment of this section covers the cases of non-multiple
components $C^{(i)}=L(p_{kl})$, $L(z)$, $L'$, or $L''$. In the next
sections we consider nonreduced components $C^{(i)}$, however the
case of multiplicity one is also covered there with the same (up to
a suitable coordinate change) answer.
\end{remark}

\subsubsection{Deformation patterns for nonisolated singularities,
I}\label{sec542} In this section we construct deformation patterns
for a component $s'L'$ of $C$ with $s'>1$.

\smallskip

{\bf Step 1: Preparation}. To relax the notation, within this
section we write $s$ for $s'$.

Consider the blow-down $\Pi:\Sig\to\PP^2$ which contracts
$E_1,...,E_6$. It takes $E$ to a conic $\Pi(E)$ passing through $5$
fixed points $\Pi(E_1),...,\Pi(E_5)$, and takes $L',L''$ to straight
lines tangent to $\Pi(E)$ and passing through the fixed point
$\Pi(E_6)$. The linear system $|D|$ on $\Sig$ turns into the linear
system $\Pi_*|D|$ of plane curves of degree $d$ with multiple points
$\Pi(E_1),...,\Pi(E_6)$ of order $d_1,...,d_6$, respectively. Take
an affine plane $\C^2\subset\PP^2$ with coordinates $x,y$ such that
$\Pi(L')=\{y=0\}$, $\Pi(E)=\{S(x,y):=y+yx+x^2=0\}$, $\Pi(E_6)$ is
the infinite point of $\Pi(L')$, and the tangency point of $\Pi(L')$
and $\Pi(E)$ is the origin. The curves $\Pi(C_t)$, where $C_t\in V$,
are then described by formula
\begin{equation}F_t(x,y)=S(x,y)\widetilde
G_t(x,y)+t^\mu G_t(x,y)\ ,\label{e31}\end{equation} which is is a
conversion of formula (\ref{e24}) (cf. also (\ref{e23})) and in
which we can suppose that $G_0(0,0)=1$ and that $\widetilde
G_0(x,y)=y^s\widetilde G'_0(x,y)$ with a polynomial $\widetilde
G'_0(x,y)$ defining the union of the components of $\Pi(C)$
different from $\Pi(L'),\Pi(L'')$. Furthermore, $a=\widetilde
G'_0(0,0)$ can be computed from relations (\ref{e22}) and
(\ref{e31}) as
\begin{equation}
a=-\frac{d}{d\tau}\left(S(x,y)
\big|_{\Pi(L_p)}\right)\Big|_{\Pi(p)}\ .\label{e32}\end{equation}
Furthermore, moving all the terms with exponent of $t$ greater or
equal to $\mu$ from the first summand in (\ref{e31}) to the second
one, we can assume that
\begin{equation}\text{all exponents of}\ t\ \text{occurring in}\
\widetilde G_t(x,y)\ \text{are strictly less than}\ \mu\
.\label{e34}\end{equation}

\smallskip

{\bf Step 2: Tropical limit}. Now we find the tropical limit of the
family (\ref{e31}) in the sense of \cite{Sh0}. For the reader's
convenience, we shortly recall what is the tropical limit. Consider
the family of the curves $K_t$ defined by polynomials (\ref{e31})
for $t\ne0$ in the trivial family of the toric surfaces
$\C(\Delta)$, where $\Delta$ is the Newton polygon of a generic
polynomial (\ref{e31}). Then we close up this family at $t=0$ (all
details can be found in \cite{Sh0}):
\begin{itemize}\item to each monomial $x^iy^j$ we assign the point
$(i,j,\bn_{ij})\in\Z^3$, where $\bn_{ij}$ is the minimal exponent of
$t$ in the coefficient of $x^iy^j$ in $F_t$, and then define a
convex piece-wise linear function $\nu:\Delta\to\R$, whose graph is
the lower part of the convex hull of all the points
$(i,j,\nu_{ij})$;
\item the surface $\C(\Delta)$ degenerates into the union of toric
surfaces $\C(\del)$, where $\del$ runs over the maximal linearity
domains of $\nu$ (the subdivision of $\Delta$ into these linearity
domains, which all are convex lattice polygons, we call
$\nu$-subdivision),\;
\item then write
$F_t(x,y)=\sum_{(i,j)\in\Delta}(a_{ij}+O(t))t^{\nu(i,j)}$ and define
the limit of the curves $K_t$ at $t=0$ as the union of the limit
curves $K_\del\subset\C(\del)$ given by limit polynomials
$f_\del=\sum_{(i,j)\in\del}a_{ij}x^iy^j$ for all pieces $\del$ of
the $\nu$-subdivision of $\Delta$.
\end{itemize} Notice that \begin{equation}\nu_{ij}\ge\nu(i,j)\
\text{for all}\ i,j,\quad \text{and}\ \nu_{ij}=\nu(i,j)\ \text{for
the vertices of}\ \nu\text{-subdivision}\ .\label{e37}\end{equation}

\smallskip

\begin{figure}
\setlength{\unitlength}{1cm}
\begin{picture}(13,20)(-2,0)
\thinlines\put(1.5,1){\vector(1,0){4}}\put(1.5,1){\vector(0,1){4.5}}
\put(7,1){\vector(1,0){4}}\put(7,1){\vector(0,1){4.5}}
\put(1.5,7.5){\vector(1,0){4}}\put(1.5,7.5){\vector(0,1){4.5}}
\put(7,7.5){\vector(1,0){4}}\put(7,7.5){\vector(0,1){4.5}}
\put(1.5,14){\vector(1,0){4}}\put(1.5,14){\vector(0,1){4.5}}
\put(7,14){\vector(1,0){4}}\put(7,14){\vector(0,1){4.5}}

\put(4.5,1){\line(0,1){1}}\put(4.5,2){\line(-1,1){3}}
\put(10,1){\line(0,1){1}}\put(10,2){\line(-1,1){3}}
\put(4.5,7.5){\line(0,1){1}}\put(4.5,8.5){\line(-1,1){3}}
\put(10,7.5){\line(0,1){1}}\put(10,8.5){\line(-1,1){3}}
\put(4.5,14){\line(0,1){1}}\put(4.5,15){\line(-1,1){3}}
\put(10,14){\line(0,1){1}}\put(10,15){\line(-1,1){3}}
\put(8,14){\line(0,1){1}}\put(1.5,7.5){\line(1,1){1}}
\put(2,8){\line(1,0){2}}\put(4,7.5){\line(0,1){0.5}}
\put(4,8){\line(1,1){0.5}}\put(7,7.5){\line(2,1){1}}
\put(7,8.5){\line(2,-1){1}}\put(2.5,1){\line(0,1){1}}
\put(2.5,1.5){\line(1,0){2}}\put(7,1.5){\line(2,-1){1}}
\put(8,1){\line(0,1){1}}

\dashline{0.2}(1.5,15)(2.5,15)\dashline{0.2}(2.5,14)(2.5,15)
\dashline{0.2}(8,7.5)(8,8.5)\dashline{0.2}(7,8)(8,8)

\thicklines \put(1.5,2.5){\line(2,-1){1}}\put(2.5,2){\line(1,0){2}}
\put(7,2.5){\line(2,-1){1}}\put(8,2){\line(1,0){2}}
\put(1.5,9){\line(2,-1){1}}\put(2.5,8.5){\line(1,0){2}}
\put(7,9){\line(2,-1){1}}\put(8,8.5){\line(1,0){2}}
\put(1.5,15.5){\line(2,-1){1}}\put(2.5,15){\line(1,0){2}}
\put(7,15.5){\line(2,-1){1}}\put(8,15){\line(1,0){2}}

\put(3.3,0){\rm (e)}\put(8.8,0){\rm (f)} \put(3.3,6.5){\rm
(c)}\put(8.8,6.5){\rm (d)}\put(3.3,13){\rm (a)}\put(8.8,13){\rm (b)}
\put(0.4,15.5){$s+1$}\put(1.1,14.9){$s$}\put(5.9,15.5){$s+1$}
\put(5.9,2.5){$s+1$}\put(6.6,1.4){$1$}\put(5.9,8.4){$j+1$}
\put(6.6,7.9){$j$}\put(2.4,13.5){$2$}\put(7.9,13.5){$2$}
\put(7.9,7){$2$}\put(2.4,0.5){$2$}\put(7.9,0.5){$2$}\put(1.1,17.9){$d$}
\put(6.6,17.9){$d$}\put(1.1,11.4){$d$} \put(6.6,11.4){$d$}
\put(1.1,4.9){$d$} \put(6.6,4.9){$d$}\put(4.2,0.5){$d-d_6$}
\put(9.7,0.5){$d-d_6$}\put(4.2,7){$d-d_6$} \put(9.7,7){$d-d_6$}
\put(4.2,13.5){$d-d_6$} \put(9.7,13.5){$d-d_6$}

\end{picture}
\caption{Deformation patterns for $sL'$}\label{f2}
\end{figure}
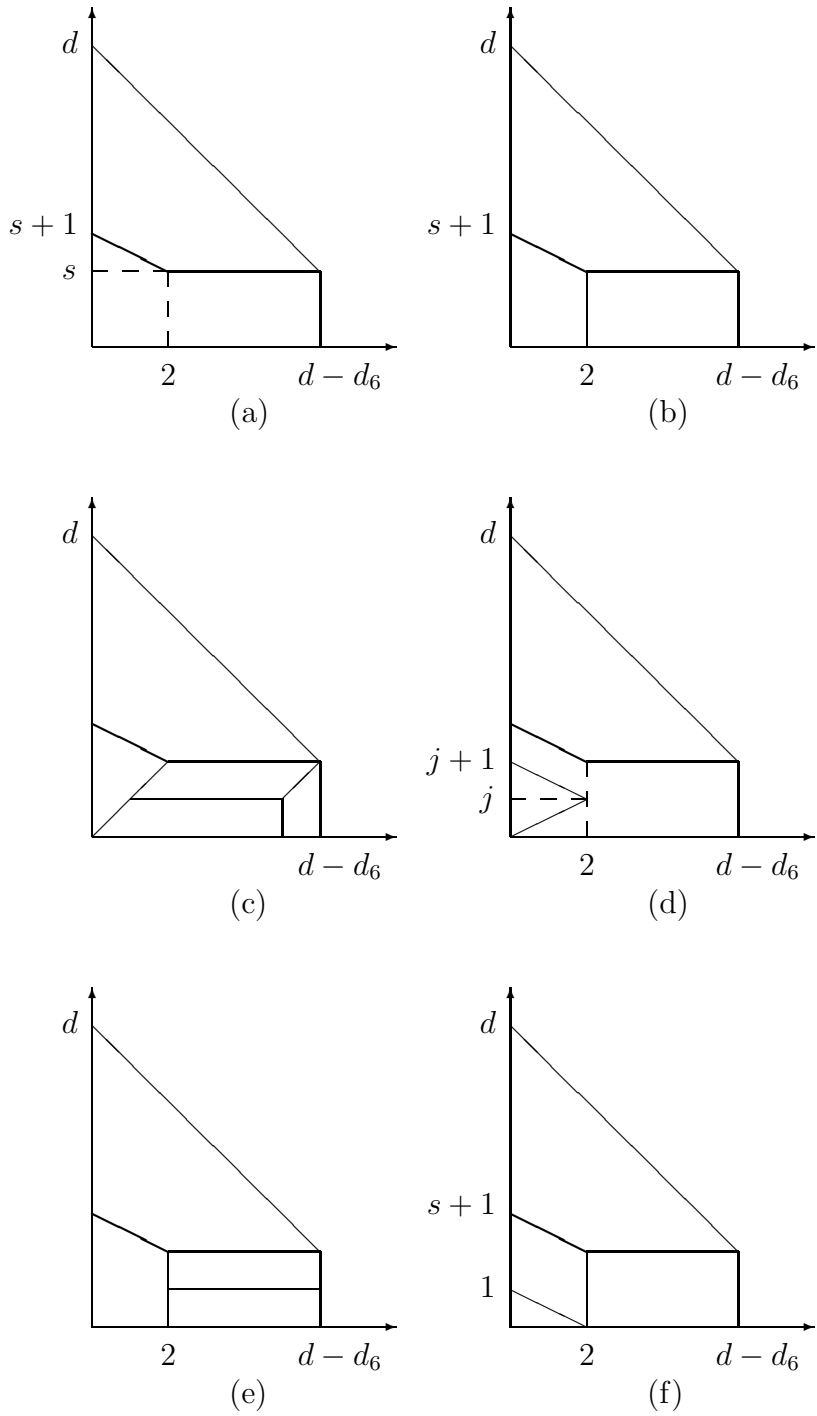

{\bf Step 3: Subdivision of the Newton polygon}. In Figure
\ref{f2}(a), we depicted the Newton polygon $\Delta$ of $F_t(x,y)$,
whose part $\del_1$ above the bold line is the Newton polygon of
$(y+yx+x^2)y^s\widetilde G'_0(x,y)$, and it is the linearity domain
of $\nu$, where it vanishes. Below the bold line, $\nu$ is positive,
and hence this part $\del_2$ is subdivided by other linearity
domains of $\nu$. We claim that $\del_2$ is subdivided into two
pieces as shown in Figure \ref{f2}(b).

{\bf (1)} The truncation of the polynomial $(y+yx+x^2)y^s\widetilde
G'_0(x,y)$ to the incline segment $\sigma_1$ of the bold line is
$(y+x^2)y^s$, and the truncation to the horizontal segment
$\sigma_2$ of the bold line is $x^2y^s\widetilde G'_0(x,0)$. Observe
that the roots of $\widetilde G'_0(x,0)$ are the coordinates of the
intersection points of the closure of $\Pi(C\setminus(L'\cup E))$
with $\Pi(L')$ (counting multiplicities). By Lemma \ref{l1}(2i), the
latter intersection points persist in the deformation $C_t$,
$t\in(\C,0)$, which implies (cf. \cite[Section 3.6]{Sh0}) that the
subdivision of $\del_2$ must contain one or few trapezes $\del$ with
a pair of horizontal sides parallel to $\sigma_2$ of length at least
$\deg\widetilde G'_0(x,0)$, and of the total height $s$ (like in
Figure \ref{f2}(c)), and the corresponding polynomials $f_\del$ are
divisible by $\widetilde G'_0(x,0)$.

{\bf (2)} Let us show that all such trapezes must be rectangles.

First, notice that $\nu(0,0)=\mu$, since
$F_t(0,0)=t^\mu(aG_0(0,0)+O(t))$. The convexity of $\nu$ and its
constancy along the horizontal sides of the trapezes (cf. Figure
\ref{f2}(c)) yield that $0\le\nu\le\mu(s-l)/s$ along such a side on
the height $0\le l<s$. Ordering the monomials of $\widetilde
G_t(x,y)$, first, by the growing power of $x$ and then by the
growing power of $y$, and using (\ref{e31}), we inductively derive
that the minimal exponents $\mu_{il}$ in the coefficients of
$x^iy^l$ in $\widetilde G_t(x,t)$ satisfy
\begin{equation}\mu_{il}\ge\nu(2,l)\quad\text{for all}\quad i\ge0,\
0\le l\le s\ .\label{e36}\end{equation} This excludes possible
vertices $v=(0,j)$ or $(1,j)$ with $0<j<s$ of the trapezes, since,
otherwise, by formula (\ref{e31}) and relation (\ref{e37}), we would
get $\nu(v)=\nu_v\ge\min\{\mu_{0,j-1},\mu_{1,j-1}\}$, contradicting
(\ref{e36}) and the strict decrease of $\nu$ with respect to the
second variable.

Next we observe that the points $(0,0)$ and $(1,0)$ cannot serve as
vertices of the lowest trapeze. Indeed, the values of $\mu_{00}$ and
$\mu_{10}$ come from the second summand in (\ref{e31}), and hence
are equal to $\mu$. Thus, if the lowest trapeze is not a rectangle,
its lower side has vertex $(0,0)$. So, suppose that $(0,0)$ is a
vertex of the lower trapeze and bring this assumption to
contradiction. The upper left vertex of that lowest trapeze must be
$(2,j)$ with some $1\le j\le s$. Thus, we have
$$0\le\nu_{2j}=\nu(2,j)\le\mu\frac{s-j}{s}\ ,$$ and by linearity of
$\nu$ in each trapeze,
$$\nu(i,l)=\nu_{2j}\frac{l}{j}+\mu\frac{j-l}{j}\quad\text{for all}\quad
i\ge 2,\ 0\le l\le j\ .$$ Combining this bound with (\ref{e36}), we
derive that
$$\nu_{il}\ge\nu_{2j}\frac{l-1}{j}+\mu\frac{j-l+1}{l}>\nu_{2j}\frac{l}{j+1}+\mu\frac{j+1-l}{j+1},\quad\text{for all}\ i=0,1,\ 1\le l\le
j\ .$$ On the other hand, from (\ref{e31}), one obtains
$\nu(0,j+1)=\mu_{0j}=\nu_{2j}$, and hence the preceding bound
together with the convexity of $\nu$ will imply that
$$\nu_{il}>\nu(i,l)\quad\text{for all}\quad i=0,1,\ 0<l\le j\ .$$
In particular, we obtain that the subdivision of $\del_2$ contains
the triangle $\quad$ $\quad$ $\quad$ $\quad$ $\quad$
\mbox{$\del=\conv\{(0,0),(2,j),(0,j+1)\}$} (see Figure \ref{f2}(d)),
and, moreover, the limit polynomial $f_\del$ contains only three
monomials corresponding to the vertices of $\del$. It is easy to see
that then the limit curve $K_\del$ is of geometric genus
$\#(\Int(\del)\cap\Z^2)>0$, which implies that the curve $C_t$ has
handles in a neighborhood of $L'$ contradicting Lemma \ref{l1}(2vi).

\smallskip

{\bf (3)} Let us show that there is only one rectangle with
horizontal/vertical edges in the subdivision of $\del_2$. Suppose
that there are several rectangles like this (cf. Figure
\ref{f2}(e)). Let $(2,j)$, $0<j<s$, be the left lower vertex of the
upper rectangle. Arguing as in the preceding paragraph, one can
derive that
\begin{eqnarray}&\nu_{1l}>\nu(1,l)\quad\text{for all}\quad 0\le l\le
s\ ,\label{e38}\\
&\nu_{0,j+1}=\nu(0,j+1)=\nu_{2,j}=\nu(2,j)=\mu_{0,j}\ ,\nonumber\\
&\nu_{0l}\ge\nu(0,l)=\nu(2,l-1)\quad\text{for all}\quad j<l\le s\
.\nonumber\end{eqnarray} From this we conclude that the
parallelogram $\del_0=\conv\{(0,j+1),(0,s+1),(2,j),(2,s)\}$ is a
part of the $\nu$-subdivision, and that the coefficients of
$f_{\del_0}$ are nonzero only along the vertical sides of $\del_0$
and they all come from the coefficients of $t^{\nu(0,l)}$ in the
coefficients of $y^{l-1}$ in $\widetilde G_t(x,y)$, $j<l\le s$.
Particularly, $f_{\del_0}=(\sum_{i=0}^{s-j}a_iy^i)(y+x^2)$, and
hence the curve $K_{\del_0}$ splits off the conic $S_0=\{y+x^2=0\}$
and $s-j$ horizontal lines, each crossing $S_0$ transversally at two
distinct points different from the origin. By a suitable variable
change $(x,y)\mapsto (xt^a,yt^b)$ in $F_t(x,y)$ and division by the
minimal power of $t$ in the obtained polynomial (still denoted by
$F_t$), we can make the corresponding function $\nu$ vanishing along
$\del_0$, which means that $F_t(x,t)=f_{\del_0}+O(t)$. At the same
time this operation takes the conic $S$ into a family of conics
$S_t=\{y+x^2+xy\cdot O(t)=0\}$. Notice that each line of
$K_{\del_0}$ crosses $S_t$ at two points convergent to its
intersection points with $S_0$. The statement of Lemma \ref{l1}(2vi)
yields that the curves $C_t$, $t\ne0$, do not intersect $E$ in a
neighborhood of $L'$. Hence, the intersection points of the lines of
$K_{\del_0}$ with $S_0$ must smooth out in the deformation
$\bigcup_\del K_\del\to K_t$, $t\ne0$. However such a smoothing will
develop at least one handle of $C_t$ in a neighborhood of $L'$,
which is a contradiction to Lemma \ref{l1}(2vi).

\smallskip

{\bf (4)} To complete the proof of the claim that the
$\nu$-subdivision of $\Delta$ is the one shown in Figure
\ref{f2}(b), it remains to exclude the subdivision depicted in
Figure \ref{f2}(f), and this can be done in the same manner as in
the preceding paragraph, where we excluded from $\nu$-subdivision
parallelograms with vertices $(0,j+1),(0,s+1),(2,j),(2,s)$.

\smallskip

{\bf Step 4: Limit curves}. As we observed above, the limit curve
corresponding to the rectangle splits into the union of vertical and
horizontal lines. For the limit polynomial $f_\del$ of the trapeze
$\del=\conv\{(0,0),(0,s+1),(2,0),(2,s)\}$, the consideration of part
(3) in Step 3, notably, (\ref{e38}), gives that $f_\del(x,y)$ does
not contain monomials $xy^l$, and the coefficients of the monomials
$y^{l+1}$ and $x^2y^l$ equal the coefficient of $t^{\nu(0,l+1)}$
standing at the monomial $y^l$ in $\widetilde G_t(x,y)$ for all
$l=0,...,s-1$. Then
\begin{equation}f_\del(x,y)=a+(y+x^2)f(y)\ ,\label{e40}\end{equation}
where $a$ is defined by (\ref{e32}), and
$f(y)=y^s+b_{s-1}y^{s-1}+...+b_0$. We call this polynomial the
\emph{deformation pattern} for the component $L'$ of $C$ associated
with the parameterized branch $V$ of $V(\overline\bp,C)$.

\begin{lemma}\label{l3}
(1) Given a curve $C=E\cup\widetilde C$ and a set $\bz'\subset
\widetilde C\cap E$ as chosen in section \ref{sec441}, the
coefficient $b$ is determined uniquely.

(2) For any fixed $b\ne0$, the set ${\mathcal P}(L')$ of the
polynomials given by (\ref{e40}) and defining a rational curve
consists of $s+1$ elements, and they all can be found from the
relation
\begin{equation}yf(y)+a=
\frac{a}{2}\left(\cheb_{s+1}\left(\frac{y}{(2^{s-1}a)^{1/(s+1)}}+y'\right)+1\right)
\ ,\label{e41}\end{equation} where $\cheb_{s+1}(y)=\cos((s+1)\arccos
y)$, and $y'$ is the only positive simple root of $\quad$
\mbox{$\cheb_{s+1}(y)-1$}.

(3) The germ at $f_\del\in{\mathcal P}(L')$ of the set of those
polynomials with Newton quadrangle
$\conv\{(0,0),(0,s+1),(2,0),(2,s)\}$ which define a rational curve,
is smooth of codimension $s$, and intersects transversally (at one
point $\{f_\del\}$) with the space of polynomials
$$\left\{(y+x^2)\left(y^s+\sum_{l=0}^{s-1}c_ly^l\right)\ :\
c_0,...,c_{s-1}\in\C\right\}\ .$$
\end{lemma}

{\bf Proof}. The first claim follows from the construction.

For the second claim, notice that $f_\del(-x,y)=f_\del(x,y)$ (cf.
(\ref{e40})); hence the curve $K_\del$ is the double cover of a
nonsingular curve ramified at the intersection points with the toric
divisors corresponding to the vertical sides of $\del$. It is easy
to show that the double cover is rational (cf. Lemma \ref{l1}(2vi)),
if and only if $f(y)$ has precisely $[s/2]$ double roots, and
$yf(y)+a$ has precisely $[(s+1)/2]$ double roots. Hence up to linear
change in the source and in the target, all such polynomials
$yf(y)+a$ must coincide with the \emph{Chebyshev polynomial}
$\cheb_{s+1}(y)$, which is characterized by the following
properties:

- it is real, has degree $s+1$, is real, and its leading coefficient
is $2^s$,

- for even $s$, it is odd, has $s/2$ maxima on the level $1$ and
$s/2$ minima on the level $-1$.

- for odd $s$, it is even, has $(s-1)/2$ maxima on the level $1$ and
$(s+1)/2$ minima on the level $-1$.

Thus, for $yf(y)+a$ we obtain precisely $s+1$ possibilities given by
formula (\ref{e41}). The smoothness and the dimension statement in
the third claim follow, for instance, from \cite[Theorem
6.1(iii)]{GK}. For the transversality statement, we notice that the
tangent space to $M$ consists of polynomials vanishing at each of
the $s$ nodes of $f_\del=0$, whereas the tangent space to the affine
space is $\{(y+x^2)\sum_{l=0}^{s-1}c_ly^l\}$, and its nontrivial
elements cannot vanish at $s$ points outside $y+x^2=0$ with distinct
$y$-coordinates.   \proofend

\begin{remark}\label{r6}
Similarly to Remark \ref{r5}(1), we observe that the values of $\nu$
must be integer at $(0,l)$, $0\le l\le s+1$, and hence $\mu$ must be
divisible by $s+1$.
\end{remark}

{\bf Proof of Proposition \ref{pdeg1}(2iv)}. In the notation of
Proposition \ref{pdeg1}, let $\bn:C^{(i)}\to L'$ is an $s$-multiple
cover, $s\ge2$. By Proposition \ref{pnov2}(ii) and Proposition
\ref{pdeg1}(2iii), $C^{(i)}$ is rational, and the cover
$\bn:C^{(i)}\to L'$ has two ramification points, one of which $z\in
L'$ differs from the tangency point of $L'$ and $E$. In particular,
if we take a tubular neighborhood $U$ of $L'$ in $\Sigma$ and the
projections $C_t\cap U\to L'$, $t\ne0$, defined by a pencil of lines
through a generic point outside $U$, then these projections will
have ramification points converging to $z$. On the other hand, it is
easy to check that the critical points of the natural projection of
the curve $K_\del\subset\C(\del)$ given by $f_\del(x,y)=0$ with
$f_\del$ as in Lemma \ref{l3} onto the toric divisor corresponding
to the lower edge of $\del$, all lie in the big torus
$(\C^*)^2\subset\C(\del)$, and hence the critical points of the
projection $C_t\cap U\to L'$ have coordinates
$((\xi+O(t))t^a,(\eta+O(t))t^{2a})$ with $a>0$, $\xi,\eta\in\C$,
thus, for $t\to0$ they converge to the origin, which in our setting
is the tangency point of $E$ and $L'$. Therefore, all ramification
points of the projection converge to that tangency point, a
contradiction. \proofend

\subsubsection{Deformation patterns for nonisolated singularities, II}\label{sec444} In this section we construct deformation patterns
for a component $C^{(i)}=kL(p_{kj})$ of $C$ with $k>1$, and we use
the technique developed  section \ref{sec542}.

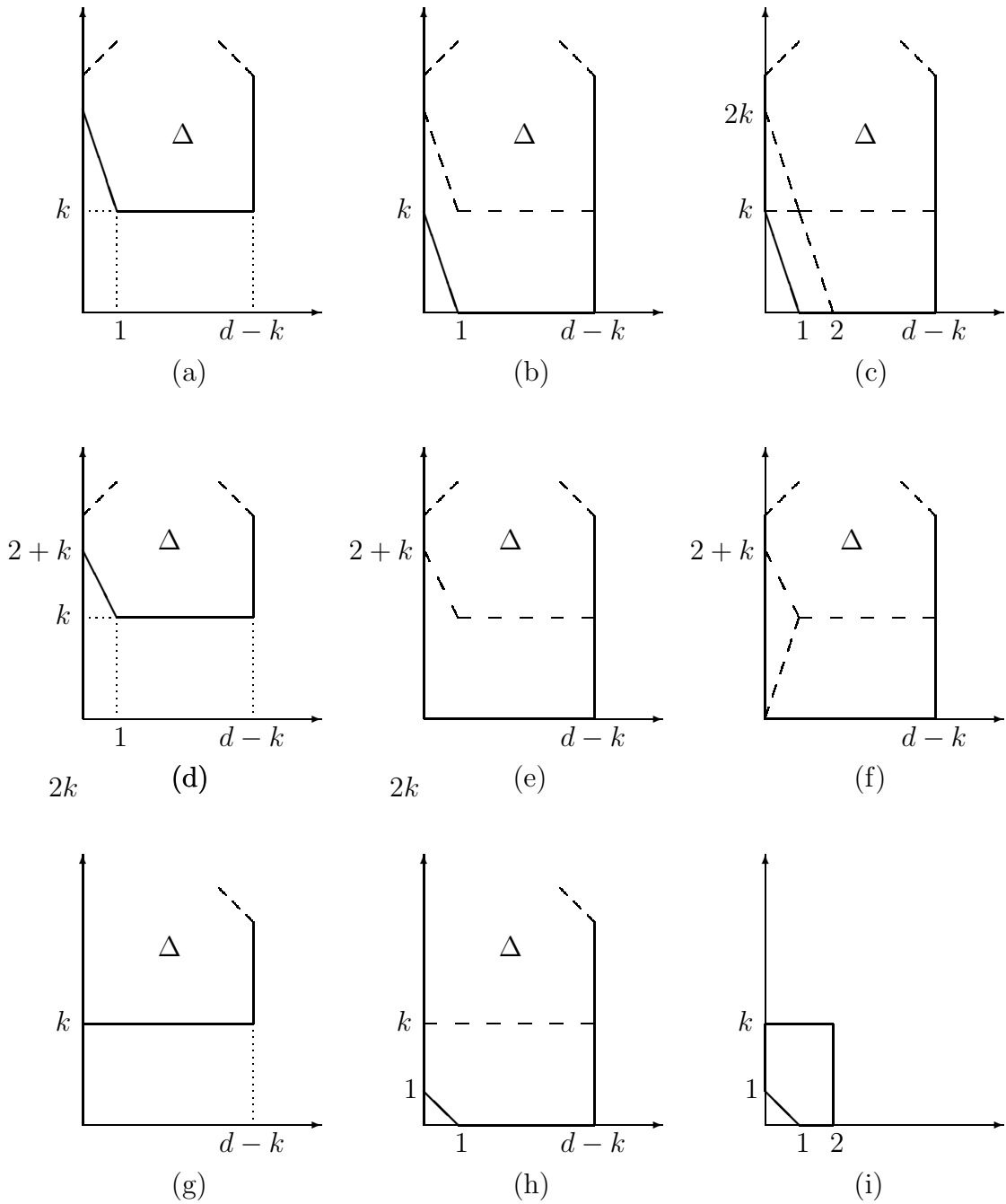
\begin{figure}
\setlength{\unitlength}{1cm}
\begin{picture}(15,18)(0,0)
\thinlines\put(1.5,1.5){\vector(1,0){3.5}}\put(6.5,1.5){\vector(1,0){3.5}}\put(11.5,1.5){\vector(1,0){3.5}}
\put(1.5,1.5){\vector(0,1){4}}\put(6.5,1.5){\vector(0,1){4}}\put(11.5,1.5){\vector(0,1){4}}
\put(1.5,7.5){\vector(1,0){3.5}}\put(6.5,7.5){\vector(1,0){3.5}}\put(11.5,7.5){\vector(1,0){3.5}}
\put(1.5,13.5){\vector(1,0){3.5}}\put(6.5,13.5){\vector(1,0){3.5}}\put(11.5,13.5){\vector(1,0){3.5}}
\put(1.5,7.5){\vector(0,1){4}}\put(6.5,7.5){\vector(0,1){4}}\put(11.5,7.5){\vector(0,1){4}}
\put(1.5,13.5){\vector(0,1){4.5}}\put(6.5,13.5){\vector(0,1){4.5}}\put(11.5,13.5){\vector(0,1){4.5}}
\thicklines\put(1.5,3){\line(0,1){1.5}} \put(1.5,3){\line(1,0){2.5}}
\put(4,3){\line(0,1){1.5}}\put(6.5,2){\line(0,1){2.5}}\put(7,1.5){\line(1,0){2}}\put(7,1.5){\line(-1,1){0.5}}
\put(9,1.5){\line(0,1){3}}\put(11.5,2){\line(0,1){1}}\put(11.5,2){\line(1,-1){0,5}}\put(12,1.5){\line(1,0){0.5}}
\put(12.5,1.5){\line(0,1){1.5}}\put(11.5,3){\line(1,0){1}}
\put(1.5,10){\line(0,1){0.5}}\put(1.5,10){\line(1,-2){0.5}}\put(2,9){\line(1,0){2}}
\put(4,9){\line(0,1){1.5}}\put(6.5,7.5){\line(0,1){3}}\put(6.5,7.5){\line(1,0){2.5}}
\put(9,7.5){\line(0,1){3}}\put(11.5,7.5){\line(0,1){3}}\put(11.5,7.5){\line(1,0){2.5}}
\put(14,7.5){\line(0,1){3}}
\put(1.5,16.5){\line(0,1){0.5}}\put(1.5,16.5){\line(1,-3){0.5}}
\put(2,15){\line(1,0){2}}\put(4,15){\line(0,1){2}}\put(6.5,15){\line(0,1){2}}
\put(6.5,15){\line(1,-3){0.5}}\put(7,13.5){\line(1,0){2}}\put(9,13.5){\line(0,1){3.5}}
\put(11.5,15){\line(0,1){2}}\put(11.5,15){\line(1,-3){0.5}}\put(12,13.5){\line(1,0){2}}
\put(14,13.5){\line(0,1){3.5}} \thinlines
\dashline{0.2}(4,4.5)(3.5,5)\dashline{0.2}(9,4.5)(8.5,5)
\dashline{0.2}(6.5,3)(9,3)
\dashline{0.2}(1.5,10.5)(2,11)\dashline{0.2}(6.5,10.5)(7,11)\dashline{0.2}(11.5,10.5)(12,11)
\dashline{0.2}(1.5,17)(2,17.5)\dashline{0.2}(6.5,17)(7,17.5)\dashline{0.2}(11.5,17)(12,17.5)
\dashline{0.2}(4,10.5)(3.5,11)\dashline{0.2}(9,10.5)(8.5,11)\dashline{0.2}(14,10.5)(13.5,11)
\dashline{0.2}(4,17)(3.5,17.5)\dashline{0.2}(9,17)(8.5,17.5)\dashline{0.2}(14,17)(13.5,17.5)
\dashline{0.2}(11.5,15)(12,15)\dashline{0.2}(12,15)(12.5,13.5)\dashline{0.2}(6.5,16.5)(7,15)
\dashline{0.2}(7,15)(9,15)\dashline{0.2}(11.5,16.5)(12,15)\dashline{0.2}(12,15)(14,15)
\dashline{0.2}(6.5,10)(7,9)\dashline{0.2}(7,9)(9,9)\dashline{0.2}(11.5,10)(12,9)
\dashline{0.2}(11.5,7.5)(12,9)\dashline{0.2}(12,9)(14,9)
\dottedline{0.1}(4,3)(4,1.5)
\dottedline{0.1}(1.5,9)(2,9)\dottedline{0.1}(2,9)(2,7.5)\dottedline{0.1}(4,9)(4,7.5)
\dottedline{0.1}(1.5,15)(2,15)\dottedline{0.1}(2,15)(2,13.5)
\dottedline{0.1}(4,13.5)(4,15) \put(2.8,6.5){\rm
(d)}\put(2.8,12.5){\rm (a)} \put(7.8,6.5){\rm (e)}\put(7.8,12.5){\rm
(b)}\put(12.8,6.5){\rm (f)}\put(12.8,12.5){\rm (c)}\put(2.8,6.5){\rm
(d)}\put(2.8,0.5){\rm (g)}\put(7.8,0.5){\rm (h)}\put(12.8,0.5){\rm
(i)}\put(2.6,10){$\Delta$}\put(2.6,4){$\Delta$}\put(7.6,4){$\Delta$}\put(7.6,10){$\Delta$}
\put(12.6,10){$\Delta$}\put(2.8,6.5){\rm
(d)}\put(2.8,16){$\Delta$}\put(7.8,16){$\Delta$}
\put(12.8,16){$\Delta$}\put(1.95,7.1){$1$}
\put(1.95,13.1){$1$}\put(6.95,1.1){$1$}\put(11.95,1.1){$1$}\put(12.45,1.1){$2$}\put(6.95,13.1){$1$}
\put(11.95,13.1){$1$}\put(12.45,13.1){$2$}
\put(1.1,8.9){$k$}\put(1.1,2.9){$k$}\put(6.1,2.9){$k$}\put(11.1,2.9){$k$}\put(6.2,1.9){$1$}\put(11.2,1.9){$1$}
\put(1.1,14.9){$k$}\put(6.1,14.9){$k$}\put(11.1,14.9){$k$}
\put(0.4,9.85){$2+k$}\put(5.4,9.85){$2+k$}\put(10.4,9.85){$2+k$}
\put(1,06.35){$2k$}\put(6,06.35){$2k$}\put(10.9,16.3){$2k$}
\put(3.5,7.1){$d-k$}\put(3.5,1.1){$d-k$}\put(8.5,1.1){$d-k$}\put(8.5,7.1){$d-k$}\put(13.5,7.1){$d-k$}
\put(3.5,13.1){$d-k$}\put(8.5,13.1){$d-k$}\put(13.5,13.1){$d-k$}
\end{picture}
\caption{Deformation patterns for $kL(p_{kl})$ and
$kL(z)$}\label{f3}
\end{figure}

Perform the blow-down $\Pi:\Sig\to\PP^2$, contracting $E_1,...,E_6$,
and take an affine plane $\C^2\subset\PP^2$ with the coordinates
$(x,y)$ such that $\Pi(L(p_{kl}))=\{y=0\}$, $\Pi(E_6)$ is the
infinitely far point of $\Pi(L(p_{kl}))$, and the conic $\Pi(E)$
crosses $\Pi(L(p_{kl}))$ transversally at $\Pi(p_{kl})=(0,0)$ and at
$(1,0)$, which is smoothed out in the deformation of $C$ along $V$.
There exists a unique transformation of $\C^2$
$$(x,\ y)\mapsto\left(x+\sum_{1\le i\le k}a_iy^i,\ y\right)\ ,$$ such that the
equation of $\Pi(E)\cap\C^2$ becomes
$$S(x,y):=y^k+x-x^2+\sum_{j>k}a_{0j}y^j+\sum_{j>0}(a_{1j}xy^j+a_{2j}x^2y^j)=0\
.$$ Thus, the curve $\Pi(C)\cap\C^2$ is defined by a polynomial of
the form $F_0(x,y)=S(x,y)\widetilde G_0(x,y)$ with Newton polygon
$\Delta$, whose lower part is the union of the segments
$[(0,2k),(1,k)]$ and $[(1,k),(d-k,k)]$ (shown by fat line in Figure
\ref{f3}(a)). Here the root of the truncation of $F_0$ on the
segment $[(0,2k),(1,k)]$ corresponds to the branch of $\Pi(E)$
centered at the origin, and the roots of the truncation of $F_0$ on
the segment $[(1,k),(d-k,k)]$ correspond to the intersection points
of $\Pi(E)$ with $\Pi(L(p_{kl}))$ (equal to $(1,0)$) and the
intersection points of $\Pi(L(p_{kl}))$ with the other components of
$\Pi(C)$. As in Step 1, section \ref{sec542}, the curves
$\Pi(C_t)\cap\C^2$ are defined by polynomials $F_t(x,y)$ converging
to $F_0(x,y)$ as $t\to0$ and given by (\ref{e31}) with the Newton
polygon $\Delta'$ being the convex hull of the union of $\Delta$
with the points $(0,k)$, $(1,0)]$, and $(d-k,0)$ (see Figure
\ref{f3}(b)).

Then, using Lemma \ref{l1}(i,iii,v) and proceeding as in Step 3,
section \ref{sec542}, we can derive that the tropical limit of the
family $F_t(x,y)$ is as follows: \begin{itemize}\item the area
$\Delta'\setminus\Delta$ is subdivided as depicted in Figure
\ref{f3}(c),
\item the limit curve $K_{\del_1}$ with the Newton parallelogram $\del_1$
splits into two components, {\it i.e.} its defining polynomial is
$(x-a')(y^s-b'x)$ with some $a',b'\ne0$, \item the limit curve
$K_{\del_2}$ with the Newton trapeze $\del_2$ is defined by the
polynomial $\quad$ \mbox{$((x-1)y^s+b'x)\psi(x)$}, where the roots
of $\psi(x)$ correspond to the intersection points (counting
multiplicities) of $\Pi(L(p_{kl}))$ with the other components of
$\Pi(C)$.
\end{itemize} The convex piece-wise linear function
$\nu:\Delta'\to\R$ defining the above subdivision (cf. Step 2,
section \ref{sec542}) is characterized by its values
$$\nu\big|_\Delta=0,\quad\nu(0,k)=\xi>0,\quad\nu(2,0)=\eta>0\ .$$ Particularly,
this yields that $\nu(1,0)=\xi+\eta$ and $\nu(u,0)=\eta$ as $u\ge
2$. Combining this with (\ref{e31}) and the fact that the curves
$C_t$, $t\ne0$, cross $E$ at $p_{kl}$ and in a neighborhood of some
intersection points of $E$ with the components of $C$ different from
$L(p_{kl})$, we derive that $\xi=\eta=\mu$ with $\mu$ defined by
(\ref{e22}), and that the limit curves are as follows:
\begin{enumerate}\item[(i)] $\{y^k-x+a'=0\}$ for the triangle
$\conv\{(0,k),(0,2k),(1,k)\}$,
\item[(ii)] $\{(x-a')(y^k-b'x)=0\}$ for the parallelogram
$\del_1=\conv\{(0,k),(1,0),(2,0),(1,k)\}$,
\item[(iii)] $\{((x-1)y^k+b'x)\psi(x)=0\}$ for the trapeze $\del_2=\conv\{(2,0),(1,k),(d-k,0),(d-k,k)\}$. \end{enumerate}

\begin{lemma}\label{l49}
The parameters $a'$ and $b'$ in the limit curves in (i)-(iii) are
uniquely defined by the curve $C$, the choice of the set $\bz'$ (cf.
section \ref{sec441}) and by relation (\ref{e22}).
\end{lemma}

{\bf Proof.} Coming back to the coordinates $(u,v)$ introduces in
section \ref{sec541}, we get that $a'$ and $b'$ linearly enter
equation (\ref{e21}) and relation (\ref{e22}). \proofend

\subsubsection{Deformation patterns for nonisolated singularities, III}\label{sec445} In this section we construct
deformation patterns for a component $C^{(i)}=kL(z)$ of $C$ with
$k>1$ and $\gam^{(i)}=e_k$, and again we apply the technique of
sections \ref{sec542}, \ref{sec444}.

Perform the blow down $\Pi:\Sig\to\PP^2$, contracting $E_1,...,E_6$,
and take an affine plane $\C^2\subset \PP^2$ with coordinates $x,y$
such that $\Pi(L(z))\cap\C^2=\{y=0\}$, $\Pi(E_6)$ is the infinitely
far point of $\Pi(L(z))$. Furthermore, we assume that the
intersection point of $\Pi(E)$ and $\Pi(L(z))$ belonging to
$\Pi(\bz')$ has coordinates $(0,0)$, the other intersection point of
$\Pi(E)$ and $\Pi(L(z))$ has coordinates $(1,0)$, and the point
$\Pi(z)\in\Pi(L_j)$ has coordinates $(x_0,0)$, $x_0\ne0,1$. At last
$\Pi(E)$ is a conic, whose equation can be made $y^2-x+x^2=0$ when
keeping the preceding data.

The curve $\Pi(C)\cap\C^2$ is then defined by a polynomial of the
form $F_0(x,y)=S(x,y)\widetilde G_0(x,y)$ with Newton polygon
$\Delta$, whose lower part is the union of the segments
$[(0,2+k),(1,k)]$ and $[(1,k),(d-k,k)]$ (shown by fat line in Figure
\ref{f3}(d)). Here the root of the truncation of $F_0$ on the
segment $[(0,2+k),(1,k)]$ corresponds to the branch of $\Pi(E)$
centered at the origin. In turn, the truncation of $F_0$ on the
segment $[(1,k),(d-k,k)]$ equals $y^kx(x-1)f(x)$, where the roots of
$f(x)$ correspond to the intersection points of $\Pi(L(z))$ with the
other components of $\Pi(C\setminus E)$. As in Step 1, section
\ref{sec542}, the curves $\Pi(C_t)\cap\C^2$ are defined by
polynomials $F_t(x,y)$ converging to $F_0(x,y)$ as $t\to0$ and given
by (\ref{e31}) with the Newton polygon
$\Delta'=\conv(\Delta\cup\{(0,0),(d-k,0)\})$ (see Figure
\ref{f3}(e)).

Then, using Lemma \ref{l1}(2i,iv,v) and proceeding as in Step 3,
section \ref{sec542}, we can derive that the tropical limit of the
family $F_t(x,y)$ is as follows: \begin{enumerate}\item[(i)] the
area $\Delta'\setminus\Delta$ is subdivided as depicted in Figure
\ref{f3}(f),
\item[(ii)] the limit curve $K_{\del_1}$ with the Newton trapeze
\begin{equation}\del_1=\conv\{(0,0),(d-k,0),(1,k),(d-k,k)\}\label{edel1}\end{equation} is defined by a
polynomial $F_{\del_1}(x,y)=(x-1)f(x)h^{(1)}_z(x,y)$ (with $f(x)$
introduced above), $h^{(1)}_z(x,y)$ has Newton triangle
$\conv\{(0,0),(1,0),(1,k)\}$, the curve $\{h^{(1)}_z(x,y)=0\}$
crosses $\Pi(L(z))$ at $\Pi(z)=(x_0,0)$ and crosses the line $x=1$
at one point with multiplicity $k$;
\item[(iii)] the convex piece-wise linear function
$\nu:\Delta\cup\del_1\cup\del_2\to\R$ whose graph projects to the
above subdivision (cf. Step 2, section \ref{sec542}) is
characterized by its values
$$\nu\big|_\Delta=0,\quad\nu(u,0)=\mu,\ 0\le u\le d-k\ .$$
\end{enumerate} Moreover, using (\ref{e34}), we can easily derive
that \begin{enumerate}\item[(iii)] the limit curve $K_{\del_2}$ with
the Newton triangle
\begin{equation}\del_2=\conv\{(0,0),(0,2+k),(1,k)\}\label{edel2}\end{equation} is defined by a
trinomial.\end{enumerate}

\begin{lemma}\label{llj}
(1) The limit curve $K_{\del_2}$, the polynomial $f(x)$ and the
coefficients $a_{00}^{z}$, $a_{10}^{z}$, $a_{1k}^{z}$ of $1$, $x$,
and $xy^k$ in $h^{(1)}_z(x,y)$ are uniquely defined by the curve
$C$, by the choice of the set $\bz'$ (cf. section \ref{sec541}), and
by relation (\ref{e22}). Within these data, there are precisely $k$
polynomials $h^{(1)}_z(x,y)$ meeting conditions (ii) above.

(2) In the space of polynomials with Newton trapeze $\del_1$ (cf.
(\ref{edel1})), the germ at $F_{\del_1}$ of the family of
polynomials, splitting into a product of binomials and a polynomial
with Newton triangle $\del_2$ (cf. (\ref{edel2})), which meets
condition (ii) above, is smooth, and it intersect with the affine
space of polynomials, having the same coefficients on the segments
$[(1,k),(d-k,k)]$ and $[(0,0),(1,k)]$ as in $F_{\del_1}$ and
vanishing at the point $(x_0,0)$, transversally at one element
(equal to $F_{\del_1}$).
\end{lemma}

{\bf Proof.} Straightforward. \proofend

We denote the set of polynomials $h^{(1)}_z$ as in Lemma \ref{llj}
as ${\mathcal P}^{(1)}(L(z))$.

\begin{remark}\label{r7}
Observe that none of the coefficients of the polynomials
$h^{(1)}_z(x,y)$ in Lemma \ref{llj} vanishes, and hence the values
$\nu(1,i)=\mu(k-i)/k$, $0\le i\le k$, must be integral, in
particular, $k$ divides $\mu$.
\end{remark}

\subsubsection{Deformation patterns for nonisolated singularities,
IV}\label{sec446} In this section, using the above approach, we
construct deformation patterns for a component $C^{(i)}=kL(z)$ of
$C$ with $k>1$ and $\gam^{(i)}=2e_k$.

Perform the blow down $\Pi:\Sig\to\PP^2$, contracting $E_1,...,E_6$,
and take an affine plane $\C^2\subset \PP^2$ with coordinates $x,y$
such that $\Pi(L(z))\cap\C^2=\{y=0\}$, $\Pi(z)=(0,0)$, $\Pi(E_6)$ is
the infinitely far point of $\Pi(L(z))$, $\Pi(L')$ is the infinitely
far line, and its tangency point with $\Pi(E)$ is the infinitely far
point of the axis $\{x=0\}$. Furthermore, the equation of $\Pi(E)$
can be brought to the form $S(x,y):=y^2+P(x)=0$ with a nondegenerate
quadratic polynomial $P(x)$. The curve $\Pi(C)\cap\C^2$ is then
defined by a polynomial of the form $F_0(x,y)=S(x,y)\widetilde
G_0(x,y)$ with Newton polygon $\Delta$, whose lower part is the
segment $[(0,k),(d-k,k)]$ (see Figure \ref{f3}(g)). Here the root of
the truncation of $F_0$ on the segment $[(0,k),(d-k,k)]$ equals
$y^kP(x)f(x)$, where the roots of $f(x)$ correspond to the
intersection points of $\Pi(L(z))$ with the other components of
$\Pi(C\setminus E)$. As in Step 1, section \ref{sec542}, the curves
$\Pi(C_t)\cap\C^2$ are defined by polynomials $F_t(x,y)$ converging
to $F_0(x,y)$ as $t\to0$ and given by (\ref{e31}) with the Newton
polygon $\Delta'=\conv(\Delta\cup\{(1,0),(0,1),(d-k,0)\})$ (see
Figure \ref{f3}(h)). Using Lemma \ref{l1}(2i,v) and proceeding as in
Step 3, section \ref{sec542}, we can derive that the tropical limit
of the family $F_t(x,y)$ is as follows:
\begin{enumerate}\item[(i)] $\Delta'$ is
subdivided as depicted in Figure \ref{f3}(h);
\item[(ii)] the limit curve $K_{\del'}$ with the Newton polygon
$\del'=\Delta'\setminus\Delta$ is defined by a polynomial
$F_{\del'}(x,y)=f(x)h_z^{(2)}(x,y)$ ($f(x)$ introduced above),
$h^{(2)}_z(x,y)$ has Newton polygon
$\conv\{(0,1),(1,0),(2,0),(2,k),(0,k)\}$;
\item[(iii)] the convex piece-wise linear function
$\nu:\Delta'\to\R$, whose graph projects to the subdivision in
Figure \ref{f3}(h), is characterized by its values
$$\nu\big|_\Delta=0,\quad\nu(u,0)=\mu,\ 1\le u\le d-k\ .$$
\end{enumerate} Using (\ref{e34}), we can easily derive that
\begin{equation}h^{(2)}_z(x,y)=y\varphi(y)P(x)+Q(x),\quad
\varphi(y)=y^{k-1}+\sum_{j=0}^{k-2}c_jy^j\ .\label{eG}\end{equation}
Finally, we observe that by Lemma \ref{l1}(2v),
$\{h^{(2)}_z(x,y)=0\}$ is a rational curve.

\begin{lemma}\label{llz}
(1) The coefficients of the polynomial $Q(x)=a_1x+a_2x^2$ in
(\ref{eG}) are uniquely defined by the curve $C$, by the choice of
the set $\bz'$ (cf. section \ref{sec541}), and by relation
(\ref{e22}). Within these data, there are precisely $k^2$
polynomials $h^{(2)}_z(x,y)$ satisfying (\ref{eG}) and defining a
rational curve.

(2) In the space of polynomials with Newton polygon $\del'$ (see
Figure \ref{f3}(h)), the germ at $F_{\del'}$ of the family of
polynomials, splitting into a product of binomials and a polynomial
with Newton polygon $\del$, which defines a rational curve, is
smooth, and it intersects transversally at one element (equal to
$F_{\del'}$) with the affine space of polynomials $g(x,y)f(x)$,
where $g(x,y)$ is given by formula (\ref{eG}) with free parameters
$c_0,...,c_{k-2}$.
\end{lemma}

{\bf Proof.} We explain only the number $k^2$ of polynomials
$h^{(2)}_z(x,y)$ in count. Computing singularities of
$\{h^{(2)}_z(x,y)=0\}$ via (\ref{eG}) and
$$h^{(2)}_z(x,y)=\frac{\partial}{\partial x}h^{(2)}_z(x,y)=\frac{\partial}{\partial y}h^{(2)}_z(x,y)=0\ ,$$
we, first, obtain an equation $QP'-Q'P=0$, which has two distinct
roots $x_1,x_2$ (it follows from the fact that, by construction, the
root of $Q(x)$ can be made generic with respect to the roots of
$P(x)$), and, second, we derive that the polynomial $y\varphi(y)$ in
(\ref{eG}) has precisely two critical levels (equal to
$-Q(x_1)/P(x_1)$ and $-Q(x_2)/P(x_2)$), and hence (cf. the proof of
Lemma \ref{l3}), up to a uniquely defined affine modification in the
source and in the target, and up to a shift, the polynomial
$y\varphi(y)$ coincides with one of the $k$ Chebyshev polynomials of
degree $k$. At last, the vanishing at $y=0$ determines precisely $k$
possible shifts for each of the above modified Chebyshev
polynomials. \proofend

We denote the set of polynomials $h^{(2)}_z$ as in Lemma \ref{llj}
as ${\mathcal P}^{(2)}(L(z))$.

\begin{remark}\label{r8}
Observe that none of the coefficients of the polynomials
$h^{(2)}_z(x,y)$ in Lemma \ref{llj} vanishes, and hence the values
$\nu(1,i)=\mu(k-i)/k$, $0\le i\le k$, must be integral, in
particular, $k$ divides $\mu$.
\end{remark}

\subsubsection{From deformation patterns to deformation} Observe that formula (\ref{e18}) counts the number of possible
collections of deformation patterns as described in Lemmas \ref{l2},
\ref{l3}, \ref{l49}, \ref{llj}, and \ref{llz} for all choices made
in section \ref{sec441}. Thus we complete the proof of Proposition
\ref{p6} with the following

\begin{lemma}\label{l4}
In the notations and hypotheses of Proposition \ref{p6}, let us be
given the set $\bz'$, chosen as in section \ref{sec441}. Then

(1) for any collection of polynomials
\begin{equation}\begin{cases}&\Psi_{ij}(x,y)\in{\mathcal
P}_i(q'_{ij})\quad \text{for all}\quad q'_{ij}\in\bz',\ \text{as in
Lemma \ref{l2}},\\ &\Psi'(x,y)\in{\mathcal P}(L')\ \text{and}\
\Psi''(x,y)\in{\mathcal P}(L''),\
\text{as in Lemma \ref{l3}},\\
&h_z^{(1)}\in{\mathcal P}^{(1)}(L(z))\quad\text{for all}\
C^{(i)}=kL(z)\ \text{with}\ \gam^{(i)}=e_k,\ \text{as in Lemma
\ref{llj}},\\ &h_z^{(2)}\in{\mathcal P}^{(2)}(L(z))\quad\text{for
all}\ C^{(i)}=kL(z)\ \text{with}\ \gam^{(i)}=2e_k,\ \text{as in
Lemma \ref{llz}},\end{cases}\label{e42}\end{equation} there exists a
unique normally parameterized branch $V$ of $V(\overline\bp,C)$
realizing the given polynomials as deformation patterns;

(2) there is one-to-one correspondence between the curves $C'\in
V(\overline\bp,C)$ passing through a given point $p'\in L_p$ and the
collections of polynomials (\ref{e42}).
\end{lemma}

{\bf Proof}. By technical reason, we shall work with
$H^0(\Sig,{\mathcal O}_\Sig(D))$ rather than with the linear system
$|D|$, and the curves $C'\in|D|$ will be lifted to sections in
$H^0(\Sig,{\mathcal O}_\Sig(D))$ denoted by $F_{C'}$ and specified
below.

Using the data of Lemma \ref{l4} and given deformation patterns
(\ref{e42}), we will derive a provisional formula for the lift of
$V$ to $H^0(\Sig,{\mathcal O}_\Sig(D))$. Then, using certain
transversality conditions, we will show that the formula does define
a unique parameterized branch $V$ of $V(\overline\bp,C)$.

\smallskip

{\bf Step 1}. In the linear system $|D|$, consider the germ at $C$
of the family of curves, which \begin{itemize}\item pass through
$\overline\bp\setminus\{z\ :\ L(z)\subset C\}$,
\item in neighborhood of the points of
$$\bp\cup(\Sing(C_{\redu})\setminus
E)\cup(\bz\setminus\{q_{kl}\in\bz\ :\ q_{kl}\in L(z),\ L(z)\subset
C\})\ ,$$ realize local deformations as described in Lemma
\ref{l1}(2i,ii,iii,iv),
\item in neighborhood of the points $\{q_{kl}\in\bz\ :\ q_{kl}\in L(z),\ L(z)\subset C\}$, realize
deformations in which the local branches of $E$ do not glue up with
local branches of the lines $L(z)$ (weaker that that in Lemma
\ref{l1}(iv), where we additionally require a special tangency
condition). \end{itemize} We will show that this family lifts to a
smooth variety germ $M\subset H^0(\Sig,{\mathcal O}_\Sig(D))$.
Namely, we will describe $M$ as an intersection of certain smooth
germs and then verify the transversality of the intersection. In the
following computations we use the model of section \ref{sec541} with
affine coordinates $u,v$.

In a neighborhood of a point
$p_{kl}\in\bp\setminus\bigcup_iC^{(i)}$, in the coordinates
$x=u-u(p_{kl})$, $y=v$, we have $$F_{C}(x,y)=y(c+\psi(x,y)),\quad
c\ne0,\ \psi(0,0)=0\ .$$ By Lemma \ref{l1}(ii), the local
deformation of $C$ along $V(\overline\bp,C)$ can be described as
$$y(c+\psi(x,y))+\sum_{kk'+l'\ge
k}c_{k'l'}x^{k'}y^{l'},\quad c_{k'l'}\in(\C,0)\ ,$$ which defines a
germ of a linear subvariety $M_{p_{kl}}\subset H^0(\Sig,{\mathcal
O}_\Sig(D))$ which embeds into ${\mathcal O}_{\Sig,p_{kl}}$ as the
intersection of the ideal $I_{p_{kl}}=\langle y,x^k\rangle$ with the
image of $H^0(\Sig,{\mathcal O}_\Sig(D))$ \footnote{There is no
canonical map $H^0(\Sig,{\mathcal O}_\Sig(D))\to{\mathcal
O}_{\Sig,p_{kl}}$, but the image of $H^0(\Sig,{\mathcal O}_\Sig(D))$
in ${\mathcal O}_{\Sig,p_{kl}}$ is defined correctly.}.

In a neighborhood of a point $p_{kl}\in\bp\cap C^{(i)}$, where
$C^{(i)}$ is reduced, in the coordinates $x=u-u(p_{kl})$, $y=v$, we
have
\begin{equation}F_{C}(x,y)=y(cy+c'x^k+\psi(x,y)),\quad
cc'\ne0,\ \psi(x,y)=\sum_{k'+kl'>2k}a_{k'l'}x^{k'}y^{l'}\
.\label{enn1}\end{equation} By Lemma \ref{l1}(ii), the local
deformation of $C$ along $V(\overline\bp,C)$ can be described as
$$\left(y+\sum_{k',l'\ge
0}c_{k'l'}x^{k'}y^{l'}\right)\left(cy+c'x^k+\psi(x,y)+\sum_{k'+kl'\ge
2k}c'_{k'l'}x^{k'}y^{l'}\right)\ ,$$ $c_{k'l'},c'_{k'l'}\in(\C,0)$.
It is immediate that this formula defines the germ of a smooth
subvariety $M_{p_{kl}}\subset H^0(\Sig,{\mathcal O}_\Sig(D))$, whose
tangent space embeds into ${\mathcal O}_{\Sig,p_{kl}}$ as the
intersection of the ideal \mbox{$I_{p_{kl}}=\langle
cy+c'x^k+\psi(x,y),y^2,x^ky\rangle$} with the image of
$H^0(\Sig,{\mathcal O}_\Sig(D))$. In the same way, for a point
$p_{kl}\in\bp\cap C^{(i)}$, where $C^{(i)}=kL(p_{kl})$, the
considered deformation as in Lemma \ref{l1}(iii) is realized in a
smooth variety germ $M_{p_{kl}}\subset H^0(\Sig,{\mathcal
O}_\Sig(D))$, whose tangent space embeds into ${\mathcal
O}_{\Sig,p_{kl}}$ as the intersection of the ideal
\mbox{$I_{p_{kl}}=\langle x^k+\psi(x,y),y^2,x^ky\rangle$} with the
image of $H^0(\Sig,{\mathcal O}_\Sig(D))$.

In a neighborhood of a point $q_{kl}\in\bz\cap C^{(i)}$, where
$C^{(i)}$ is reduced, in the coordinates $x=u-u(q_{kl})$, $y=v$, we
again have (\ref{enn1}), and by Lemma \ref{l1}(iv), the local
deformation of $C$ along $V(\overline\bp,C)$ can be described as
$$\left(y+\sum_{k',l'\ge
0}c_{k'l'}x^{k'}y^{l'}\right)\left(cy+d(x+c')^k+\psi(x,y)+\sum_{k'+kl'\ge
2k}c'_{k'l'}x^{k'}y^{l'}\right)\ ,$$
$c_{k'l'},c',c'_{k'l'}\in(\C,0)$. This formula also defines the germ
of a smooth subvariety $\quad\quad$ \mbox{$M_{q_{kl}}\subset
H^0(\Sig,{\mathcal O}_\Sig(D))$}, whose tangent space embeds into
${\mathcal O}_{\Sig,q_{kl}}$ as the intersection of the ideal
\mbox{$I_{q_{kl}}=\langle cy+dx^k+\psi(x,y),y^2,x^{k-1}y\rangle$}
with the image of $H^0(\Sig,{\mathcal O}_\Sig(D))$.

For a point $q_{kl}\in\bz\cap C^{(i)}$, where $C^{(i)}=kL(z)$, in
the local coordinates $x=u-u(q_{kl})$, $y=v$, where, in addition,
$L(z)=\{x=0\}$, the considered deformation can be described as
$$\left(y+\sum_{k',l'\ge
0}c_{k'l'}x^{k'}y^{l'}\right)\left(x^k+\sum_{k',l'\ge
0}c'_{k'l'}x^{k'}y^{l'}\right)\ ,$$ $c_{k'l'},c'_{k'l'}\in(\C,0)$.
This formula again defines the germ of a smooth subvariety
$\qquad\qquad$ \mbox{$M_{q_{kl}}\subset H^0(\Sig,{\mathcal
O}_\Sig(D))$}, whose tangent space embeds into ${\mathcal
O}_{\Sig,q_{kl}}$ as the intersection of the ideal
\mbox{$I_{q_{kl}}=\langle y,x^k\rangle$} with the image of
$H^0(\Sig,{\mathcal O}_\Sig(D))$.

At last, we notice that, if $z\in\Sig\setminus E$ is a center of two
transverse smooth local branches $P_1,P_2$ of $C$ of multiplicities
$r_1,r_2$, respectively, and, in some local coordinates $x,y$, we
have $z=(0,0)$, $P_i=\{f_i(x,y)=0\}$, $i=1,2$, then (cf.
\cite[Section 5.2]{Sh0}) the closure $M_z\subset H^0(\Sig,{\mathcal
O}_\Sig(D))$ of the set of sections defining in a neighborhood of
$z$ curves with $r_1r_2$ nodes, is smooth at $F_{C}$ and its tangent
space embeds into ${\mathcal O}_{\Sig,z}$ as the intersection of the
ideal $I_z=\langle f_1^{r_1}(x,y),f_2^{r_2}(x,y)\rangle$ with the
image of $H^0(\Sig,{\mathcal O}_\Sig(D))$.

The branch $V$ we are looking for must lie inside the intersection
of the above germs $M_{p_{kl}},M_{q_{kl}},M_z$ with the space of
sections vanishing at $\overline\bp\setminus\{z\ :\ L(z)\subset
C\}$. We claim that this intersection is transverse, or,
equivalently, that $H^1(\Sig,{\mathcal J}_{Z_1/\Sig}(D))=0$, where
${\mathcal J}_{Z_1/\Sig}\subset{\mathcal O}_\Sig$ is the ideal sheaf
of the scheme $Z_1\subset\Sig$ supported at
\mbox{$\bp\cup\bz\cup(\Sing(C_{\redu})\setminus E)$}, where it is
defined by the above local ideals $I_{p_{kl}},I_{q_{kl}},I_z$,
respectively, and supported at $\overline\bp\setminus\{z\ :\
L(z)\subset C\}$, where it is defined by the maximal ideals. From
the exact sequence of sheaves (cf. \cite{Hir})
$$0\to{\mathcal J}_{Z_1:E/\Sig}(D-E)\to{\mathcal J}_{Z_1/\Sig}(D)\to{\mathcal J}_{Z_1\cap E/E}(DE)\to0$$
we obtain the cohomology exact sequence
\begin{equation}H^1(\Sig,{\mathcal J}_{Z_1:E/\Sig}(D-E))\to H^1(\Sig,{\mathcal J}_{Z_1/\Sig}(D))\to H^1(E,{\mathcal J}_{Z_1\cap
E/E}(DE))\ .\label{enn3}\end{equation} Here $H^1(E,{\mathcal
J}_{Z_1\cap E/E}(D))=0$ in view of Riemann-Roch, since $\deg(Z_1\cap
E)=DE$. Thus, it remains to prove
\begin{equation}H^1(\Sig,{\mathcal J}_{Z_1:E/\Sig}(D-E))=0\
,\label{enn4}\end{equation} Notice that the scheme $Z_1:E$ coincides
with $Z_1$ in $\Sig\setminus E$, and it is defined at the points
$p_{kl}\in\bp\cap C^{(i)}$ by the ideals $J_{p_{kl}}=\langle
y,x^k\rangle\subset{\mathcal O}_{\Sig,p_{kl}}$, at the points
$q_{kl}\in\bz\setminus\{q_{kl}\ :\ q_{kl}\in L(z),\ L(z)\subset C\}$
by the ideals $J_{q_{kl}}=\langle y,x^{k-1}\rangle\subset{\mathcal
O}_{\Sig,q_{kl}}$. From the definition of the ideals $I_z$,
$z\in\Sing(\widetilde C_{\redu})\setminus E$, and the Noether
fundamental theorem, we have a canonical decomposition
$$H^0(\Sig,{\mathcal J}_{Z_1:E/\Sig}(D-E))\simeq
\bigoplus_{C^{(i)}\ne s'L',s''L''}\prod_{j\ne i}F_{C^{(j)}}\cdot
(F_{L'})^{s'}(F_{L''})^{s''}\cdot H^0(\Sig,{\mathcal
J}_{Z_1(C^{(i)})/\Sig}(C^{(i)}))$$
$$\qquad\qquad\qquad\oplus\prod_{i=1}^mF_{C^{(i)}}\cdot (F_{L''})^{s''}\cdot
H^0(\Sig,{\mathcal O}_\Sig((s')(L-E_6)))$$
\begin{equation}\qquad\qquad\qquad\;\oplus\prod_{i=1}^mF_{C^{(i)}}\cdot
(F_{L'})^{s'}\cdot H^0(\Sig,{\mathcal O}_\Sig((s'')(L-E_6)))\
,\label{enn13}\end{equation} where $Z_1(C^{(i)})\subset C^{(i)}$ is
the part of $Z_1:E$ supported along $C^{(i)}$ and disjoint from the
other components of $C$ different from $E$ (in particular,
$Z_1(kL(z))=\emptyset$). It follows that (\ref{enn4}) turns into a
sequence of equalities
\begin{equation}
H^1(\Sig,{\mathcal J}_{Z_1(C^{(i)})/\Sig}(C^{(i)}))=0,\quad
i=1,...,m\ . \label{enn12}\end{equation} We intend to prove even a
stronger result: \begin{itemize} \item for a reduced
$C^{(i)}\not\in|L-E_6|$, we replace (\ref{enn12}) with
\begin{equation}H^1(\Sig,{\mathcal
J}_{Z_2(C^{(i)})/\Sig}(C^{(i)}))=0\ ,\label{enn8}\end{equation}
where the zero-dimensional scheme $Z_2(C^{(i)})$ is obtained from
$Z_1(C^{(i)})$ by adding the (simple) points $\overline\bp\cap
C^{(i)}$ and the scheme supported at $\bz'\cap C^{(i)}$ and defined
in $q'_{kl}\in\bz'\cap C^{(i)}$ by the ideal $J_{q'_{kl}}=\langle
y,x^{k-1}\rangle\subset{\mathcal O}_{\Sig,q'_{kl}}$; \item for
$C^{(i)}=kL(z)$ such that $\gam^{(i)}=e_k$, we replace
(\ref{enn12}), which in view of $Z_1(kL(z))=\emptyset$ reads as
$H^1(\Sig,{\mathcal O}_\Sig(kL(z)))=0$, with
\begin{equation}H^1(\Sig,{\mathcal J}_{Z^1_2(kL(z))/\Sig}(kL(z)))=0\ ,\label{enn8a}\end{equation} where the zero-dimensional
scheme $Z^1_2(kL(z))$ is supported at the point $q_{kl}\in
L(z)\cap\bz$ and is defined there by the ideal $\langle
y,x^k\rangle$; \item for $C^{(i)}=kL(z)$ such that
$\gam^{(i)}=2e_k$, we replace (\ref{enn12}) with
\begin{equation}H^1(\Sig,{\mathcal J}_{Z^2_2(kL(z))/\Sig}(kL(z)))=0\
,\label{enn8b}\end{equation} where the zero-dimensional scheme
$Z^2_2(kL(z))$ is supported at the point $z\in L(z)\cap\overline\bp$
and is defined there by the ideal $\langle x,y^k\rangle$ in the
coordinates $x,y$ introduced in section \ref{sec446}.
\end{itemize}

To derive (\ref{enn8}), we notice that $H^1(\Sig,{\mathcal
J}_{Z_1(C^{(i)})\setminus\overline\bp/\Sig}(C^{(i)}))=0$ is
equivalent to \begin{equation}H^1(\hat C^{(i)},{\mathcal N}_{\hat
C^{(i)}}(-\bd^{(i)}))=0\ ,\label{enn6}\end{equation} where $\hat
C^{(i)}$ is the normalization of $C^{(i)}$, and
$$\bd^{(i)}=\sum_{p_{kl}\in\bp\cap C^{(i)}}k\cdot
\hat p_{kl}+\sum_{q_{kl}\in\bz\cap C^{(i)}}(k-1)\cdot \hat
q_{kl}+\sum_{q'_{kl}\in\bz'\cap C^{(i)}}(k-1)\cdot \hat q'_{kl}$$
($\hat q'_{kl}$ being the lift of $q'_{kl}$ to $\hat C^{(i)}$). In
its turn, (\ref{enn6}) comes from (\ref{enn5}). This yields, in
particular, that (details of computation left to the reader)
$$h^0(\Sig,{\mathcal
J}_{Z_1(C^{(i)}\setminus\overline\bp/\Sig}(C^{(i)}))=h^0(\Sig,{\mathcal
O}_\Sig(C^{(i)}))-\deg(Z_1(C^{(i)})\setminus\overline\bp)=\#(\overline\bp\cap
C^{(i)})+1\ .$$ Due to the generic choice of the configuration
$\overline\bp$, we conclude that $$h^0(\Sig,{\mathcal
J}_{Z_2(C^{(i)})/\Sig}(C^{(i)}))=h^0(\Sig,{\mathcal
J}_{Z_1(C^{(i)})\setminus\overline\bp/\Sig}(C^{(i)}))-\#(\overline\bp\cap
C^{(i)})=1\ ,$$ which immediately implies (\ref{enn8}). Similarly,
relations (\ref{enn12}) for $C^{(i)}=kL(p_{kl})$ and relations
(\ref{enn8a}), (\ref{enn8b}) for $C^{(i)}=kL(z)$ lift to the
$k$-sheeted coverings $\bn:\hat L(p_{kl})\to L(p_{kl})$, $\bn:\hat
L(z)\to L(z)$ (cf. Proposition \ref{pnov2}(iii,iv)) in the form
$$H^1(\hat L(p_{kl}),{\mathcal N}_{\hat L(p_{kl})}(-\bn^*(p_{kl}))=H^1(\hat
L(z),{\mathcal N}_{\hat L(z)}(-\bn^*(q_{kl})))=H^1(\hat
L(z),{\mathcal N}_{\hat L(z)}(-\bn^*(z)))=0\ ,$$ which all come from
Riemann-Roch due to
$$\deg\bn^*(p_{kl})=\deg\bn^*(q_{kl})=\deg\bn^*(z)=k\ .$$

Thus, we conclude that the considered family $M$ is smooth.

\smallskip

{\bf Step 2}. Now we construct a certain parametrization of $M$.
For, choose a special basis for its tangent space
$H^0(\Sig,{\mathcal J}_{Z_1/\Sig}(D))$.

\smallskip

{\bf (2a)} In view of (\ref{enn4}), we have
\begin{equation}H^0(\Sig,{\mathcal J}_{Z_1/\Sig}(D))\simeq F_E\cdot
H^0(\Sig,{\mathcal J}_{Z_1:E/\Sig}(D-E))\oplus H^0(E,{\mathcal
J}_{Z_1\cap E/E}(DE))\ ,\label{enn2606}\end{equation} and hence
there is a section $F\in H^0(\Sig,{\mathcal
J}_{Z_1\cup\overline\bp/\Sig}(D))$, which in the coordinates $u,v$
of section \ref{sec541} can be written in the form
$$f_2(u)+v\widetilde f_2(u,v)\ ,$$ where $f_2$ is given by
(\ref{enn10}).

\smallskip

{\bf (2b)} Using (\ref{enn8}) for a reduced $C^{(i)}\not\in|L-E_6|$,
we get $H^0(\Sig,{\mathcal
J}_{Z_2(C^{(i)})/\Sig}(C^{(i)}))=F_{C^{(i)}}\cdot\C$, where
$F_{C^{(i)}}$ is specified so that, in the coordinates $x,y$ in a
neighborhood of each point $q'_{kl}\in\bz'\cap C^{(i)}$ (see section
\ref{sec541}), the coefficient of $y^2$ in $F_{C^{(i)}}(x,y)$ equals
$a_{02}^{kl}$ as defined in section \ref{sec541}. Next, from the
definition of the scheme $Z_2(C^{(i)})$, we get that
$H^0(\Sig,{\mathcal J}_{Z_1(C^{(i)})/\Sig}(C^{(i)}))$, $1\le l\le
m$, admits a basis consisting of $F_{C^{(i)}}$ and the sections
$F_{i,j,kl}$ labeled by the points $q'_{kl}\in\bz'\cap C^{(i)}$ and
the numbers $0\le j<k-1$, such that \begin{itemize}\item
$\jet_{k-2,q'_{kl}}(F_{i,j,kl}(x,0))=x^j$, $0\le j<k-1$, in the
coordinates $x=u-u(q'_{kl})$, $y=v$:
\item
$\jet_{k'-2,q'_{k'l'}}(F_{i,j,kl}(x,0))=0$ in the coordinates
$x=u-u(q'_{k'l'})$, $y=v$, for all $i,j,k,l$ and
$q'_{k'l'}\in\bz'\cap C^{(i)}$, $(k',l')\ne(k,l)$.
\end{itemize}

\smallskip

{\bf (2c)} The scheme $Z_1(kL(p_{kl}))$ is determined by the
condition to cross $E$ at the point $p_{kl}$ with multiplicity $k$.
Hence $H^0(\Sig,{\mathcal
J}_{Z_1(kL(p_{kl}))/\Sig}(kL(p_{kl})))=(F_{L(p_{kl})})^k\cdot\C$.

\smallskip

{\bf (2d)} For the $k$-dimensional space $H^0(\Sig,{\mathcal
O}_\Sig(C^{(i)}))$, where $C^{(i)}=kL(z)$, we choose the basis as
follows: \begin{itemize}\item if $\gam^{(i)}=e_k$, we take
$(F_{L(z)})^j(F_{\widetilde L(z)})^{k-j}$, $0\le j\le k$, where, in
the coordinates $x,y$ of section \ref{sec445}, $\widetilde
L(z)\in|L-E_6|$ is the infinitely far line (respectively,
$F_{L(z)}(x,y)=y$, $F_{\widetilde L(z)}(x,y)=1$);
\item if $\gam^{(i)}=2e_k$, we take $(F_{L(z)})^j(F_{\widetilde L(z)})^{k-j}$, $0\le j\le k$, where, in the
coordinates $x,y$ of section \ref{sec446}, $\widetilde
L(z)\in|L-E_6|$ is the infinitely far line (respectively,
$F_{L(z)}(x,y)=y$, $F_{\widetilde L(z)}(x,y)=1$).\end{itemize}

\smallskip

{\bf (2e)} Finally, for $H^0(\Sig,{\mathcal O}_\Sig(s'(L-E_6)))$ and
$H^0(\Sig,{\mathcal O}_\Sig(s''(L-E_6)))$ (cf. (\ref{enn13})),
choose the following bases
$$F_{L'}^kF_{\widetilde L'}^{s'-j},\ 0\le j\le s',\quad\text{and}\quad
F_{L''}^kF_{\widetilde L''}^{s''-j},\ 0\le j\le s''\ ,$$ where
$F_{L'}=y$, $F_{\widetilde L'}(x,y)=1$ in the coordinates $x,y$
introduced in section \ref{sec542}, Step 1, and $F_{L''}$,
$F_{\widetilde L''}$ are defined similarly.

\smallskip

Thus, our basis for $H^0(\Sig,{\mathcal J}_{z_1/\Sig}(D))$ contains
the section $F$ (see (2a)), the sections built in (2b)-(2e)
(multiplied with extra factors along (\ref{enn13}),
(\ref{enn2606})), and certain unspecified completion.

We can write then a parametrization of $M$ in the form
$$F_E\cdot\prod_{C^{(i)}\ne
kL({pkl}),kL(z)}\left(F_{C^{(i)}}+\sum_{ij,k}\xi_{j,k,l}F_{i,j,kl}\right)\cdot\prod_{C^{(i)}=kL(p_{kl})}(F_{L(p_{kl})})^k$$
$$\times
\prod_{C^{(i)}=kL(z)}\left((F_{L(z)})^k+\sum_{j=1}^{k}\xi_{j,z}(F_{L(z)})^{k-j}(F_{\widetilde
L(z)})^j\right)$$
$$\quad\;\times\left(F_{L'}^{s'}+
\sum_{j=0}^{s'-1}\xi'_jF_{L'}^jF_{\widetilde
L'}^{s'-j}\right)\left(F_{L''}^{s''}+
\sum_{j=0}^{s''-1}\xi''_jF_{L''}^jF_{\widetilde
L''}^{s''-j}\right)$$
\begin{equation}+\xi_0F+O(\overline\xi^2)\
,\label{enn14}\end{equation} where $\overline\xi$ is the sequence of
\begin{equation}N:=\sum_{C^{(i)}\ne kL(z)}(I\gam^{(i)}-|\gam^{(i)}|)+\sum_{C^{(i)}=kL(z)}k+s'+s''+1\label{eparam}\end{equation} free parameters
$$\xi_0,\xi_{j,s},\xi_{ijkl},\xi'_s,\xi''_s\in(\C,0)\ .$$

\smallskip

{\bf Step 3}. We are ready to construct the desired parameterized
branch $V$ of $V(\overline\bp,C)$ realizing deformation patterns
(\ref{e42}). Namely, imposing additional conditions to the
parameters $\overline\xi$ in (\ref{enn14}), we will expose a formula
for a lift of the required branch $V$ to $H^0(\Sig,{\mathcal
O}(D))$.

Put
$$\mu=\lcm\left(\{k\
:\ q'_{kl}\in\bz'\}\cup\{s'+1,s''+1\}\cup\{k\ :\
kL(z)\in\{C^{(i)}\}_{i=1,...,m}\}\right)$$ (cf. Remarks \ref{r5},
\ref{r6}, \ref{r7}), and \ref{r8}) and let $\xi_0=t^\mu a(t)$ with
$t\in(\C,0)$ and $a(0)=a$ given by formula (\ref{e23}).

\smallskip

{\bf (3a)} For a point $q'_{kl}\in\bz\cap C^{(i)}$ such that
$C^{(i)}\not\in|L-E_6|$ is reduced, choose a deformation pattern
$$\Psi_{kl}(x,y)=a^{kl}_{02}y^2+a^{kl}_{j1}x^jy+\sum_{j<k,\ j\equiv
k\mod2}b^{kl}_{j1}x^jy+b^{kl}_{00}$$ (cf. (\ref{e39}), (\ref{e42}),
and Lemma \ref{l2}(2)). Let $F_{kl}(x,y)$ be the section
(\ref{enn14}) written in the coordinates $x=u-u(q'_{kl})$, $y=v$.
The coefficient $a_{k-1,1}^{kl}(\overline\xi)$ of $x^{k-1}y$ in
$F_{kl}$ vanishes at $\overline\xi=0$. Since the coefficient of
$x^ky$ in $F_{kl}$ turns into $a_{k1}^{kl}\ne0$ as $\overline\xi=0$,
there is a function $\tau_{kl}(\overline\xi)$ vanishing at the
origin such that the coefficient of $x^{k-1}y$ in
$F_{kl}(x-\tau(\overline\xi),y)$ equals $0$, and hence the
coefficient of $x^jy$, $j<k-1$, in $F_{kl}(x-\tau(\overline\xi),y)$
equals $\xi_{i,j,kl}+\lambda_{i,j,kl}(\overline\xi)$ with a function
$\lambda_{i,j,kl}\in O(\overline\xi^2)$. Thus, due to the smoothness
and transversality statement in Lemma \ref{l2}(3), the conditions to
fit the given deformation pattern $\Psi_{ij}$ and to realize a local
deformation as described in Lemma \ref{l1}(v) amounts in a system of
relations
\begin{equation}\xi_{i,j,kl}+\lambda_{i,j,kl}(\overline\xi)=\begin{cases}t^{\mu(k-j)/(2k)}(b_{j1}^{kl}+O(t)),\quad
& j\equiv k\mod2,\\ O(t^{[\mu(k-j)/(2k)]+1}),\quad & j\not\equiv
k\mod2,\end{cases},\quad 0\le j<k-1\ .\label{enn15}\end{equation}

\smallskip

{\bf (3b)} Let $C^{(i)}=kL(z)$ with $\gam^{(i)}=e_k$, and let, in
the coordinates $x,y$ of section \ref{sec445},
$h^{(1)}_z(x,y)=x(y^k+\sum_{j=0}^{k-1}b_{z,j}y^j)-x_0b_{z,0}\in{\mathcal
P}^{(1)}(L(z))$. In these coordinates, expression (\ref{enn14})
reads as
$$(y^2-x+x^2)\left(y^k+\sum_{j=0}^{k-1}\xi_{z,j}y^j\right)(1+O(x,y,\overline\xi))+t^\mu(x_0b_{z,0}+O(t,\overline\xi))(x-1+O(x^2,y,\overline\xi))\ .$$
In view of the smoothness and transversality statement in Lemma
\ref{llj}(2), the conditions to fit the deformation pattern
$h^{(1)}_z(x,y)$, to realize a local deformation as described in
Lemma \ref{l1}(2iv), and to hit the point $z\in
L(z)\cap\overline\bp$ altogether amount in a system of equations
\begin{equation}\xi_{z,j}+O(\overline\xi^2)=t^{\mu(k-j)/k}(b_{z,j}+O(t,\overline\xi)),\quad
0\le j<k\ .\label{elj}\end{equation} Similarly, let $C^{(i)}=kL(z)$
with $\gam^{(i)}=2e_k$, and let, in the coordinates of section
\ref{sec446},
$h^{(2)}_z(x,y)=P(x)(y^k+\sum_{j=1}^{k-1}c_{z,j}y^j)+Q(x)\in{\mathcal
P}^{(2)}(L(z))$ be a deformation pattern (cf. (\ref{eG})). In these
coordinates, expression (\ref{enn14}) reads as
$$(y^2+P(x))\left(y^k+\sum_{j=0}^{k-1}\xi_{z,j}y^j\right)(1+O(x,y,\overline\xi))+t^\mu(F\big|_z+O(t,x,y,\overline\xi))\
.$$ Again the smoothness and transversality statement in Lemma
\ref{llz}(2) convert the condition to fit the deformation pattern
$h^{(2)}_z(x,y)$ and to realize a deformation of the union of
$kL(z)$ with the germs of $E$ at the points $E\cap L(z)$ into in
immersed cylinder, into a system of equations
\begin{equation}\xi_{z,j}+O(\overline\xi^2)=t^{\mu(k-j)/k}(b_{z,j}+O(t,\overline\xi)),\
1\le j<k,\quad\xi_{z,0}=t^{\mu}(-F\big|_z+O(t,\overline\xi))\
.\label{elz}\end{equation}

\smallskip

{\bf (3c)} Let $\Psi'(x,y)\in{\mathcal P}(L')$ be a deformation
pattern for the component $(L')^{s'}$ of $C$. In agreement with
(\ref{e40}) we can write (slightly modifying notations)
$$\Psi'(x,y)=b'+(y+x^2)f'(y),\quad f'(y)=y^{s'}+b'_{s'-1}y^{s'-1}+...+b'_0\
,$$ with $f'(y)$ defined as in Lemma \ref{l3}. On the other hand, in
the coordinates $x,y$ introduced in section \ref{sec542}, Step 1,
expression (\ref{enn14}) with $\xi_0=t^\mu(a+O(t))$ reads as
$$(y+xy+x^2)\left(y^{s'}+\sum_{k=0}^{s'-1}\xi'_ky^k+O(\overline\xi^2)\right)(1+O(x,y,\overline\xi))+t^\mu(b'+O(t,\overline\xi))\ .$$
Again in view of the smoothness and transversality statement in
Lemma \ref{l3}(3), the conditions to fit the given deformation
pattern $\Psi'$ together with the demand to realize a deformation in
a neighborhood of the point $L'\cap E$ as described in Lemma
\ref{l1}(2vi), turn into a system of equations
\begin{equation}\xi'_k+O(\overline\xi^2)=t^{\mu(s'+1-k)/(s'+1)}(b'_k+O(t,\overline\xi)),\quad k=0,...,s'\
.\label{enn16}\end{equation} Similarly, for the component
$(L'')^{s''}$ of $C$ with a given deformation pattern $\Psi''$, we
obtain equations
\begin{equation}\xi''_k+O(\overline\xi^2)=t^{\mu(s''+1-k)/(s''+1)}(b''_k+O(t,\overline\xi)),\quad k=0,...,s''\
.\label{enn17}\end{equation}

\smallskip

{\bf (3d)} In the right-hand side of $N-1=\sum_{C^{(i)}\ne
kL(z)}(I\gam^{(i)}-|\gam^{(i)}|)+\sum_{C^{(i)}=kL(z)}k+s'+s''$ (cf.
(\ref{eparam})) equations (\ref{enn15}), (\ref{elj}), (\ref{enn16}),
and (\ref{enn17}), we have $N-1$ unknown functions. Observe that the
linearization of this unified system has a nondegenerate diagonal
form. Hence it is uniquely soluble as far as we take
$\xi_0=t^\mu(a+O(t))$. This provides a uniquely defined
parameterized branch $V$ of $V(\overline\bp,C)$ matching all the
given deformation patterns.

\smallskip

Had we taken $\xi_0=t^{r\mu}(a+O(t))$ with some $r>1$ in the
presented construction, we were coming (due to the uniqueness claim)
to the branch of $V(\overline\bp,C)$, geometrically coinciding with
$V$, with the parametrization obtained from that of $V$ via the
parameter change $t\mapsto t^r$.

The proof of Lemma \ref{l4} and Proposition \ref{p6} is completed.
\proofend

\section{Recursive formula}
\label{sec:formula}

The set $\Pic(\Sig,E) \times \Z \times \Z_+^\infty \times
\Z_+^\infty$ is a semigroup with respect to the operation
$$(D',g',\alp',\bet')+(D'',g'',\alp'',\bet'')=(D'+D'',g'+g''-1,\alp'+\alp'',\bet'+\bet'')\ .$$
Denote by $A(\Sig,E)\subset\Pic(\Sig,E) \times \Z \times \Z_+^\infty
\times \Z_+^\infty$ the subsemigroup generated by
\begin{itemize}\item all the quadruples $(D,g,\alp,\bet)$ with reduced,
irreducible $D$, $g\ge0$, and $\alp,\bet$ satisfying (\ref{en1}),
and \item the quadruples $(s(L-E_6),0,e_s,e_s)$ and
$(s(L-E_6),0,0,2e_s)$ for all $s\ge2$.
\end{itemize} Notice that one may have $g<0$ for elements
$(D,g,\alp,\bet)\in A(\Sig,E)$; in such a case we always assume
$N_\Sig(D,g,\alp,\bet)=0$.

\begin{theorem}\label{t1}
Let $\Sig = \PP^2_6$. Given an element $(D,g,\alp,\bet)\in
A(\Sig,E)$ such that
$$R_\Sig(D,g,\bet)>0\quad\text{and}\quad(D,g,\alp,\bet)\ne(s(L-E_6),0,0,2e_s),\ s>1\ ,$$ we have
$$ N_\Sig(D,g,\alp,\bet) = \sum_{j \geq 1,\ \bet_j > 0} j \cdot
N_\Sig(D,g,\alp+e_j,\bet-e_j)$$
\begin{equation}
+ \sum \binom{\alp}{\alp^{(1)},\ldots,\alp^{(m)}}
\frac{(n-1)!}{n_1!...n_m!}\binom{k+3}{3} \prod_{i\in S} \left(
\binom{\bet^{(i)}}{\gam^{(i)}} I^{\gam^{(i)}}
N_\Sig(D^{(i)},g^{(i)},\alp^{(i)},\bet^{(i)}) \right),
\label{eq:main}\end{equation} where
$$n =R_\Sig(D,g,\bet),\quad n_i = R_\Sig(D^{(i)},g^{(i)},\bet^{(i)}),\quad   i =1,...,m\
,$$ $$S=\{i\in[1,m]\ :\
(D^{(i)},g^{(i)},\alp^{(i)},\bet^{(i)})\ne(s(L-E_6),0,e_s,e_s),\
s\ge1\}\ ,$$ and the second sum is taken
\begin{itemize}
\item   over all integer $k\ge0$ and all splittings in $A(\Sig,E)$
\begin{equation} \label{eq:splittingsq}
(D-E-k(L-E_6),g',\alp',\bet') = \sum_{i=1}^m
(D^{(i)},g^{(i)},\alp^{(i)},\bet^{(i)})\ ,\end{equation} such that
$$g^{(i)}\ge0\quad\text{and}\quad(D^{(i)},g^{(i)},\alp^{(i)},\bet^{(i)})\ne(s(L-E_6),0,0,se_2),\
s\ge 1,\ i=1,...,m\ ,$$ of all possible quadruples
$(D-E-k(L-E_6),g',\alp',\bet') \in A(\Sig,E)$ with $k\ge0$,
satisfying
\begin{enumerate}
\item[(a)] $\alp' \leq \alp$, $\bet' \geq \bet$, $g - g' = \|\bet' - \bet\| + 1$,
\item[(b)] each summand $(D^{(i)},g^{(i)},\alp^{(i)},\bet^{(i)})$ with $n_i = 0$ and $\alp^{(i)}=0$ appears in
(\ref{eq:splittingsq}) at most once,
\end{enumerate}
\item over all splittings
\begin{equation} \label{eq:splittingsbeta}
\bet' = \bet + \sum_{i=1}^m \gam^{(i)},\quad\|\gam^{(i)}\|
> 0, \ 1 \leq i \leq m
\end{equation}
satisfying the restrictions $\bet^{(i)} \geq \gam^{(i)}$ for all $1
\leq i \leq m$,
\end{itemize}
and, finally, the second sum in (\ref{eq:main}) is factorized by
simultaneous permutations in both splittings (\ref{eq:splittingsq})
and (\ref{eq:splittingsbeta}).
\end{theorem}

{\bf Proof}. Immediately follows from from Propositions \ref{pdeg1},
\ref{p5}, and \ref{p6}. So, the first summand corresponds to the
degeneration described in Proposition \ref{pdeg1}(1), and its
multiplicity $j$ was computed in Proposition \ref{p5}. The second
sum corresponds to the degeneration described in Proposition
\ref{pdeg1}(2), and its multiplicity was computed in Proposition
\ref{p6}. Notice only that the subtraction of $k(L-E_6)$ in
(\ref{eq:splittingsq}) corresponds to all possible combinations
$s'L'\cup s''L''$, $s'+s''=k$, which, by Proposition \ref{p6},
altogether contribute to the multiplicity $(k+1)+k\cdot 2+(k-1)\cdot
3+...+(k+1)=\binom{k+3}{k}$. \proofend

\smallskip

Note that the second sum in the right hand side of (\ref{eq:main})
becomes empty if $D - E$ is not effective.

\begin{theorem} \label{prop:rigidcases}
Let $\Sig = \PP^2_6$, then all the numbers $N_\Sig(D,g,\alp,\bet)$
with $(D,g,\alp,\bet) \in A(\Sig,E)$ and $R_\Sig(D,g,\bet)>0$ are
determined recursively by formula (\ref{eq:main}) and the initial
values $N_\Sig(D,g,\alp,\bet)$ indicated in Proposition \ref{pini}.
\end{theorem}

{\bf Proof}. Straightforward. \proofend

\begin{example}
For the popular case of $D=6L-2(E_1+...+E_6)$ (cf. \cite[Page 119,
Table II]{DFI}, \cite[Page 88]{GP}, \cite[Section 9.2]{Va}), our
formula gives Gromov-Witten invariants for all genera:
\begin{center}
\begin{tabular}{|l||c|c|c|c|c|c|}
\hline $g$ &\makebox[0.35cm] 0 &\makebox[0.35cm] 1 &\makebox[0.35cm]
2&\makebox[0.35cm] 3 &\makebox[0.35cm] 4
  \\
\hline
$N$ & 3240 & 1740 & 369 & 33 & 1  \\
\hline
\end{tabular}
\end{center}
\end{example}

In the important case of genus zero, all components of degenerate
curves in Proposition \ref{pdeg1} are rational, and the general
formula (\ref{eq:main}) reduces to the count of only genus zero
terms. Furthermore, it can be simplified so that the terms
corresponding to multiple covers will not explicitly appear.

Let us write for short $N_\Sig(D,\alp,\bet):=N_\Sig(D,0,\alp,\bet)$
and introduce the semigroup
$A_0(\Sig,E)\subset\Pic(\Sig,E)\times(\Z_+^{\infty})^2$, generated
by the triples $(D,\alp,\bet)$ with reduced, irreducible $D$ and
$\|\alp+\bet\|=DE$.

\begin{corollary}\label{czero}
Let $\Sig = \PP^2_6$. For any element $(D,\alp,\bet)\in A_0(\Sig,E)$
such that $R_\Sig(D,0,\bet)>0$, we have $$ N_\Sig(D,\alp,\bet) =
\sum_{j \geq 1,\ \bet_j > 0} j \cdot N_\Sig(D,\alp+e_j,\bet-e_j)$$
$$ + \sum
\frac{2^{\|\bet^{(0)}\|}I^{\bet^{(0)}}}{\bet^{(0)}!}\binom{k+3}{3}\binom{\alp}{\alp^{(0)},\ldots,\alp^{(m)}}
\binom{n-1}{n_1,\ldots,n_m}$$ \begin{equation}\times \prod_{i=1}^m
\left(\binom{\bet^{(i)}}{\gam^{(i)}} I^{\gam^{(i)}}
N_\Sig(D^{(i)},\alp^{(i)},\bet^{(i)}) \right)\
,\label{eqrat}\end{equation} where
$$n =R_\Sig(D,0,\bet),\quad n_i = R_\Sig(D^{(i)},0,\bet^{(i)}),\quad   i =1,...,m\
,$$ and the second sum is taken
\begin{itemize} \item over
all integers $k\ge0$ and vectors
$\alp^{(0)},\bet^{(0)}\in\Z_+^\infty$ such that $\alp^{(0)}\le\alp$,
$\bet^{(0)}\le\beta$;
\item over all sequences
\begin{equation}(D^{(i)},\alp^{(i)},\bet^{(i)})\in A_0(\Sig,E),\
1\le i\le m \ ,\label{eq:splittings0}\end{equation} such that
\begin{enumerate}
\item[(i)] for
all $i=1,...,m$, $D^{(i)}\ne L-E_6$ and
$R_\Sig(D^{(i)},0,\bet^{(i)})\ge0$,
\item[(ii)]
$D-E=\sum_{i=1}^mD^{(i)}+(k+I\alp^{(0)}+I\bet^{(0)})(L-E_6)$,
\item[(iii)] $\sum_{i=0}^m\alp^{(i)}\le\alp$, $\sum_{i=0}^m\bet^{(i)}\ge\bet$,
\item[(iv)] each
triple $(D^{(i)},0,\bet^{(i)})$ with $n_i=0$ appears in
(\ref{eq:splittings0}) at most once,
\end{enumerate}
\item over all sequences \begin{equation}\gam^{(i)}\in\Z_+^\infty,\quad
\|\gam^{(i)}\|=1,\quad i=1,...,m\ ,\label{eq:splittings1}
\end{equation} satisfying
$$\bet^{(i)}\ge\gam^{(i)}, \ i=1,...,m,\quad\text{and}\quad
\sum_{i=1}^m\left(\bet^{(i)}- \gam^{(i)}\right)=\bet-\bet^{(0)}\ ,$$
\end{itemize}
and the second sum in (\ref{eqrat}) is factorized by simultaneous
permutations in the sequences (\ref{eq:splittings0}) and
(\ref{eq:splittings1}).
\end{corollary}

\begin{example}
We would like to illustrate the factor $\binom{k+3}{3}$ in formulas
(\ref{eq:main}) and (\ref{eqrat}), coming from the count of
deformations of multiple lines $L',L''\in|L-E_6|$, tangent to $E$.
Namely, we compute $GW_0(\PP^2_6,(n+2)L-nE_1)$ in two ways.

Using the evident equality
$$GW_0(\PP^2_6,(n+2)L-nE_1)=GM_0(\PP^2_6,(n+2)L-nE_1-E_2-E_3-E_4)$$
and a suitable Cremona transformation, we get (in the notations of
Corollary \ref{czero})
$$GW_0(\PP^2_6,(n+2)L-nE_1)=GW_0(\PP^2_6,(2n+1)L-n(E_1+...+E_4))
=N_\Sig(D_n,0,2e_1)\ ,$$ $$D_n= (2n+1)L-n(E_1+...+E_4), \ n\ge0\ .$$
Formula (\ref{eqrat}), in which we ignore $E_5,E_6$, gives
\begin{align*}N_\Sig(D_n,0,2e_1)=&N_\Sig(D_n,e_1,e_1)=N_\Sig(D_n,2e_1,0)+2N_\Sig(D_{n-1},0,2e_1),\\
N_\Sig(D_n,2e_1,0)=&2N_\Sig(D_{n-1},e_1,e_1)+2N_\Sig(D_{n-1},0,e_2)\\
&+4\sum_{i=0}^{n-1}\binom{2n-1}{2i}N_\Sig(D_i-L+E_1,0,e_1)
N_\Sig(D_{n-i-1}(D_{n-i-1}-E_1,0,e_1)\\
=&2N_\Sig(D_{n-1},0,2e_1)+2N_\Sig(D_{n-1},0,e_2)+4\sum_{i=0}^{n-1}\binom{2n-2}{2i},\\
N_\Sig(D_{n-1},0,e_2)=&2M_\Sig(D_{n-1},e_2,0)\\
=&8\sum_{i=0}^{n-2}
\binom{2n-3}{2i}N_\Sig(D_i-L+E_1,0,e_1)N_\Sig(D_{n-i-1}-E_1,0,e_1)\\
=&8\sum_{i=0}^{n-2} \binom{2n-3}{2i}\ ,\end{align*} which implies
\begin{align*}N_\Sig(D_n,0,2e_1)=&4N_\Sig(D_{n-1},0,2e_1)+4\sum_{i=0}^{n-1}\binom{2n-2}{2i}+8\sum_{i=0}^{n-2}
\binom{2n-3}{2i}\\ =&4N_\Sig(D_{n-1},0,2e_1)+4^n(n+1)\ ,\end{align*}
and hence
$$GW_0(\PP^2_6,(n+2)L-nE_1)=4^n\binom{n+2}{2}\ .$$ On the other
hand, $$GW_0(\PP^2_6,(n+2)L-nE_1)=GW_0(\PP^2_6,(n+2)L-nE_6)=
N_\Sig(D,0,(2n+4)e_1)\ ,$$ where $D=(n+2)L-nE_6$. We again apply
formula (\ref{eqrat}):
$$N_\Sig(D,0,(2n+4)e_1)=N_\Sig(D,e_1,(2n+3)e_1)=...=N_\Sig(D,(n+4)e_1,ne_1)\ ,$$ then for $m\le n$
\begin{align*}N_\Sig(D,&(2n+4-m)e_1,me_1)=N_\Sig(D,(2n+5-m)e_1,(m-1)e_1)\\
&+2^m\sum_{k=0}^{n-m}\binom{k+3}{3}\binom{2n+4-m}{n-m-k}N_\Sig(L-E_6,0,2e_1)^m
N_\Sig(L-E_6,e_1,e_1)^{n-m-k}\\
=&N_\Sig(D,(2n+5-m)e_1,(m-1)e_1)+2^m\sum_{k=0}^{n-m}\binom{k+3}{3}\binom{2n+4-m}{n-m-k}\
,\end{align*} and hence
$$GW_0(\PP^2_6,(n+2)L-nE_1)=\sum_{m=0}^n2^m\sum_{k=0}^{n-m}\binom{k+3}{3}\binom{2n+4-m}{n-m-k}\
.$$ A direct proof of the equality
$$4^n\binom{n+2}{2}=\sum_{m=0}^n2^m\sum_{k=0}^{n-m}\binom{k+3}{3}\binom{2n+4-m}{n-m-k}$$
was communicated to us by S. Lando, to whom we are very grateful for
that. We leave it to the reader as an exercise.
\end{example}

\vskip10pt

{\ncsc Institute of Mathematics, \\[-21pt]

Hebrew University of Jerusalem, \\[-21pt]

Givat Ram, Jerusalem, 91904, Israel} \\[-21pt]

\emph{E-mail}: {\ntt mendy.shoval@mail.huji.ac.il}

\vskip10pt

{\ncsc School of Mathematical Sciences, \\[-21pt]

Raymond and Beverly Sackler Faculty of Exact Sciences,\\[-21pt]

Tel Aviv University,
Ramat Aviv, 69978 Tel Aviv, Israel} \\[-21pt]

\emph{E-mail}: {\ntt shustin@post.tau.ac.il}

\end {document}